\documentclass{article}

\title{\LARGE \textbf{On the circumference, connectivity and dominating cycles}}
\author{Zh.G. Nikoghosyan\footnote{G.G. Nicoghossian (up to 1997)}\\
Institute for Informatics and Automation Problems\\ National Academy of Sciences\\
P.Sevak 1, Yerevan 0014, Armenia\\ E-mail: zhora@ipia.sci.am}

\begin{document}

\maketitle

\begin{abstract}
Every 4-connected graph with minimum degree $\delta$ and connectivity $\kappa$ either has a cycle of length at least $4\delta-2\kappa$ or has a dominating cycle.\\
  

\end{abstract}

\section{Introduction}

Long and large cycles (paths) are the main general research objects in hamiltonian graph theory basing on two different initial conceptions: Hamilton and dominating cycles (paths). A cycle $C$ of a graph $G$ is called a Hamilton cycle if it contains every vertex of $G$. A cycle $C$ is called a dominating cycle if $V(G-C)$ is an independent set.

In 1952, Dirac [2] obtained the first sufficient condition for a graph to be hamiltonian by showing that every graph $G$ of order $n$ and minimum degree $\delta$ at least $n/2$, is hamiltonian. Although this bound $n/2$ is tight in Dirac's theorem, it was essentially lowered in several ways by direct incorporation of some additional graph invariants (namely connectivity $\kappa$ and independence number $\alpha$) due to the author [6, 7], Nash-Williams [5] and the author [8]. The latter recently was improved by Yamashita [15]. The evolution of these developments can be demonstrated by the following order.\\

\noindent\textbf{Theorem A} [2]. Every graph with $\delta\geq\frac{1}{2}n$ is hamiltonian.\\

\noindent\textbf{Theorem B} [6, 7]. Every 2-connected graph with $\delta \geq \frac{1}{3}(n+\kappa)$ is hamiltonian.\\

\noindent\textbf{Theorem C} [5]. Every 2-connected graph with $\delta \geq \max\{\frac{1}{3}(n+2),\alpha\}$ is hamiltonian.\\

\noindent\textbf{Theorem D} [8]. Every 3-connected graph with $\delta\geq\max\lbrace \frac{1}{4}(n+2\kappa),\alpha\rbrace$ is hamiltonian.\\

\noindent\textbf{Theorem E} [15]. Every 3-connected graph with $\delta\geq\max\lbrace \frac{1}{4}(n+\kappa+3),\alpha\rbrace$ is hamiltonian.\\

A short proof of Theorem B was given in [3] due to H\"{a}ggkvist. Theorems A-D have their reverse versions due to Dirac [2], the author [7], Voss and Zuluaga [14] and the author [9], respectively, concerning the circumference and Hamilton cycles.\\  

\noindent\textbf{Theorem F} [2]. Every 2-connected graph has a cycle of length at least $\min\lbrace n,2\delta\rbrace$.\\

\noindent\textbf{Theorem G} [7]. Every 3-connected graph has a cycle of length at least $\min\lbrace n,3\delta -\kappa\rbrace$.\\

\noindent\textbf{Theorem H} [14]. Every 3-connected graph with $\delta \geq \alpha$ has a cycle of length at least $\min \lbrace n, 3\delta -3\rbrace$.\\

\noindent\textbf{Theorem I} [9, 11]. Every 4-connected graph with $\delta\geq\alpha$ has a cycle of length at least $\min\lbrace n,4\delta-2\kappa\rbrace$.\\

Theorem I was firstly announced in 1985 [9]. The detailed proof of this theorem recently was given in [11]. In view of Theorem E, the following reverse version is reasonable.\\

\noindent\textbf{Conjecture 1}. Every 4-connected graph with $\delta\geq\alpha$ has a cycle of length at least $\min\lbrace n,4\delta-\kappa-4\rbrace$.\\

It is natural to look for analogous results for dominating cycles. Observe that Theorems C, D, E, H and I can not have dominating versions due to the condition $\delta\geq\alpha$. In 1971, Nash-Williams [5] obtained the first Dirac-type sufficient condition for dominating cycles, i.e. the dominating version of Theorem A. The dominating version of Theorem B is due to Sun, Tian and Wei [12], which recently was improved by Yamashita [15].\\

\noindent\textbf{Theorem J} [5]. Every 2-connected graph with $\delta \geq \frac{1}{3}(n+2)$ has a dominating cycle.\\

\noindent\textbf{Theorem K} [12]. Every 3-connected graph with $\delta \geq \frac{1}{4}(n+2\kappa)$ has a dominating cycle.\\

\noindent\textbf{Theorem L} [15]. Every 3-connected graph with $\delta \geq \frac{1}{4}(n+\kappa+3)$ has a dominating cycle.\\

The reverse version of Theorem J is due to Voss and Zuluaga [14]. \\

\noindent\textbf{Theorem M} [14]. Every 3-connected graph either has a cycle of length at least $3\delta-3$ or has a dominating cycle.\\

In this paper we prove the reverse version of Theorem K, recently conjectured in [10]. \\

\noindent\textbf{Theorem 1}. Every 4-connected graph either has a cycle of length at least $4\delta-2\kappa$ or has a dominating cycle.\\

The limit example $4K_2+K_3$ shows that the 4-connectivity condition in Theorem 1 can not be replaced by 3-connectedness. Further, the limit example $5K_2+K_4$ shows that the conclusion "has a cycle of length at least $4\delta-2\kappa$" in Theorem 1 can not be replaced by "has a cycle of length at least $4\delta-2\kappa$+1". Finally, the limit example $H(1,n-2\delta,\delta,\kappa)$ shows that the conclusion "has a dominating cycle" can not be replaced by "has a Hamilton cycle", where $H(1,n-2\delta,\delta,\kappa)$ is defined as follows. Given four integers $a,b,t,\kappa$ with $\kappa\leq t$, we denote by $H(a,b,t,\kappa)$ the graph obtained from $tK_a+\overline{K}_t$ by taking any $\kappa$ vertices in subgraph $\overline{K}_t$ and joining each of them to all vertices of $K_b$.\\

In view of Theorem E and Conjecture 1, the following reverse version of Theorem L is reasonable.\\

\noindent\textbf{Conjecture 2}. Every 4-connected graph either has a cycle of length at least $4\delta-\kappa-4$ or has a dominating cycle.\\

Let $G$ be a graph. For $X\subset V(G)$, we denote by $N(X)$ the set of all vertices of $G-X$ having neighbors in $X$. Furthermore, $\hat{X}$ is defined as $V(G)-(X\cup N(X))$. Following Hamidoune [4], we define a subset $X$ of $V(G)$ to be a fragment of $G$ if $N(X)$ is a minimum cut-set and $\hat{X}\neq\emptyset$. If $X$ is a fragment then $\hat{X}$ is a fragment too and $\hat{\hat{X}}=X$. For convenience, we will use $X^{\uparrow}$ and $X^{\downarrow}$ to denote $X$ and $\hat{X}$, respectively. Throughout the paper, we suppose w.l.o.g. that $|A^{\uparrow}|\geq|A^{\downarrow}|$. An endfragment is a fragment that contains no other fragments as a proper subset. 

To prove Theorem 1, we present four more general Dirac-type results for dominating cycles centered around a lower bound $c\geq4\delta-2\kappa$ under four alternative conditions in terms of endfragments.\\

\noindent\textbf{Theorem 2}. Let $G$ be a 3-connected graph. If $|A^{\uparrow}|\leq 3\delta-\kappa-4$ and $|A^{\downarrow}|\leq 3\delta-3\kappa+1$ for an endfragment $A^{\downarrow}$ of $G$, then either $G$ has a cycle of length at least $4\delta-2\kappa$ or has a dominating cycle. \\

\noindent\textbf{Theorem 3}. Let $G$ be a 4-connected graph. If $|A^{\uparrow}|\leq 3\delta-\kappa-4$ and $|A^{\downarrow}|\geq 3\delta-3\kappa+2$ for an endfragment $A^{\downarrow}$ of $G$, then either $G$ has a cycle of length at least $4\delta-2\kappa$ or has a dominating cycle.\\

\noindent\textbf{Theorem 4}. Let $G$ be a 4-connected graph. If $|A^{\uparrow}|\geq 3\delta-\kappa-3$ and $|A^{\downarrow}|\leq 3\delta-3\kappa+1$ for an endfragment $A^{\downarrow}$ of $G$, then either $G$ has a cycle of length at least $4\delta-2\kappa$ or has a dominating cycle. \\

\noindent\textbf{Theorem 5}. Let $G$ be a 4-connected graph. If $|A^{\uparrow}|\geq 3\delta-\kappa-3$ and $|A^{\downarrow}|\geq 3\delta-3\kappa+2$ for an endfragment $A^{\downarrow}$ of $G$, then either $G$ has a cycle of length at least $4\delta-2\kappa$ or has a dominating cycle. \\

\section{Definitions and notations}

By a graph we always mean a finite undirected graph $G$ without loops or multiple edges. A good reference for any undefined terms is [1]. For $H$ a subgraph of $G$ we will denote the vertices of $H$ by $V(H)$ and the edges of $H$ by $E(H)$. For every $S\subset V(G)$ we use $G-S$ short for $\langle V(G)-S\rangle$, the subgraph of $G$ induced by $V(G)-S$. In addition, for a subgraph $H$ of $G$ we use $G-H$ short for $G-V(H)$.  

Let $\delta$ denote the minimum degree of vertices of $G$. The connectivity $\kappa$ of $G$ is the minimum number of vertices whose removal from $G$ results in a disconnected or trivial graph. We say that $G$ is $s$-connected if $\kappa\ge s$. A set $S$ of vertices is independent if no two elements of $S$ are adjacent in $G$. The cardinality of maximum set of independent vertices is called the independence number and denoted by $\alpha$. 

Paths and cycles in a graph $G$ are considered as subgraphs of $G$. If $Q$ is a path or a cycle of $G$, then the length of $Q$, denoted by $|Q|$, is $|E(Q)|$. Throughout this paper, the vertices and edges of $G$ can be interpreted as cycles of lengths 1 and 2, respectively. A graph $G$ is hamiltonian if it contains a Hamilton cycle. The length of a longest cycle (the circumference) will be denoted by $c$.

Let $Q$ be a cycle of $G$ with a fixed cyclic orientation. For $x\in V(Q)$, we denote the successor and predecessor of $x$ on $Q$ by $x^+$ and $x^-$, respectively. For $X\subseteq V(Q)$, we define $X^+=\{x^+|x\in X\}$ and $X^-=\{x^-|x\in X\}$.  For two vertices $x$ and $y$ of $Q$, let $x\overrightarrow{Q}y$ denote the segment of $Q$ from $x$ to $y$ in the chosen direction on $Q$ and $x\overleftarrow{Q}y$ denote the segment in the reverse direction. We also use similar notation for a path $P$ of $G$. We call $x\overrightarrow{Q}y$ an $m$-segment if $|x\overrightarrow{Q}y|\geq m$. For $P$ a path of $G$, denote by $F(P)$ and $L(P)$ the first and the last vertices of $P$, respectively.   

Let $Q$ be a cycle of a graph $G$, $r\geq 2$ a positive integer and $Z_1,Z_2,...,Z_p$ are subsets of $V(Q)$ with $p\geq 2$. A collection $(Z_{1},...,Z_{p})$ is called a $(Q,r)$-scheme if $x\overrightarrow{Q}y$ is a 2-segment for each distinct $x,y\in Z_{i}$ (where $i\in\lbrace1,...,p\rbrace)$ and is an $r$-segment for each distinct $x\in Z_{i}$ and $y\in Z_{j}$ (where $i,j\in\{1,...,p\}$ and $i\neq j)$. A $(Q,r)$-scheme is nontrivial if $(Z_{1},...,Z_{p})$ has a system of distinct representatives. The definition of $(Q,r)$-scheme was first introduced by Nash-Williams [5] for $p=2$. \\

\noindent\textbf{Definition A $\{Q^{\uparrow}_{1},...,Q^{\uparrow}_{m};V^{\uparrow}_{1},...,V^{\uparrow}_{m};V^{\uparrow}\}$}. Let $A^{\uparrow}$ be a fragment of $G$ with respect to a minimum cut-set $S$. Define $Q^{\uparrow}_{1},...,Q^{\uparrow}_{m}$ as a collection of vertex disjoint paths in $\langle A^{\uparrow}\cup S\rangle$ with terminal vertices in $S$ such that 

$(*1)$\quad $|V(Q^{\uparrow}_i)|\geq 2$ $(i=1,...,m)$ and $\sum_{i=1}^{m}|V(Q^{\uparrow}_{i})|$ is as great as possible. 

Abbreviate $V^{\uparrow}_{i}=V(Q^{\uparrow}_{i})$ $(i=1,...,m)$ and $V^{\uparrow}=\bigcup^{m}_{i=1}V^{\uparrow}_{i}$. \\

\noindent\textbf{Definition B $\{Q^{\downarrow}_{1},...,Q^{\downarrow}_{m};Q^{\downarrow}_{*};V^{\downarrow}_{1},...,V^{\downarrow}_{m};V^{\downarrow}\}$}. Let $A^{\uparrow}$ be a fragment of $G$ with respect to a minimum cut-set $S$ and let $Q^{\uparrow}_{1},...,Q^{\uparrow}_{m}$ be as defined in Definition A.

If $|A^{\downarrow}|\geq2(\delta-\kappa+1)$, then we denote by $Q^{\downarrow}_{1},...,Q^{\downarrow}_{m}$ a collection of paths (if exist) in $\langle A^{\downarrow}\cup S\rangle$ with 

$(*2)$ combining $Q^{\uparrow}_{1},...,Q^{\uparrow}_{m}$ with $Q^{\downarrow}_{1},...,Q^{\downarrow}_{m}$ results in a simple cycle such that $\sum_{i=1}^{m}|V(Q^{\downarrow}_{i})|$ is as great as possible. 

Abbreviate $V^{\downarrow}_{i}=V(Q^{\downarrow}_{i})$ $(i=1,...,m)$ and $V^{\downarrow}=\bigcup^{m}_{i=1}V^{\downarrow}_{i}$. For the special case $|V^{\downarrow}\cap S|=2$ and $S\not\subseteq V^{\uparrow}$, we will use $Q^{\downarrow}_{*}$ to denote a longest path in $\langle A^{\downarrow}\cup S\rangle$ such that combining $Q^{\uparrow}_{1}$ and $Q^{\downarrow}_{*}$ results in a simple cycle with $|V(Q^{\downarrow}_{*})\cap S|\geq3$.

For the reverse case  $|A^{\downarrow}|\leq2\delta-2\kappa+1$, we use the same notation $Q^{\downarrow}_{1},...,Q^{\downarrow}_{m}$ (if no ambiguity can arise) to denote a collection of paths in $\langle A^{\downarrow}\cup S\rangle$ with  

$(*3)$ combining $Q^{\uparrow}_{1},...,Q^{\uparrow}_{m}$ with $Q^{\downarrow}_{1},...,Q^{\downarrow}_{m}$ results in a simple cycle such that $|(\cup_{i=1}^{m}V^{\downarrow}_{i})\cap S|$ is as great as possible,

$(*4)$ $|\cup_{i=1}^{m}V^{\downarrow}_{i}|$ is as great as possible, subject to $(*3)$.

In this case, $V^{\downarrow}_{i}$ and $V^{\downarrow}$ are as defined in the previous case.\\

\noindent\textbf{Definition C}. Throughout the paper we will use $C$ to denote the cycle obtained by combining $Q^{\uparrow}_{1},...,Q^{\uparrow}_{m}$ and $Q^{\downarrow}_{1},...,Q^{\downarrow}_{m}$. Assume w.l.o.g. that 

$(*5)$\quad $Q^{\uparrow}_{1},...,Q^{\uparrow}_{m}$ is chosen such that $|C|$ is as great as possible.

\section{Preliminaries}

In [5], Nash-Williams proved the following result concerning $(Q,r)$-schemes for a cycle $Q$ and a pair $(Z_1 ,Z_2 )$ of subsets of $V(Q)$.\\

\noindent\textbf{Lemma A [5]}. Let $Q$ be a cycle and $(Z_{1},Z_{2})$ be a nontrivial $(Q,r)$-scheme. Then
$$
|Q|\geq\min\Big\{ 2(|Z_{1}|+|Z_{2}|)+2r-6,\frac{1}{2}r(|Z_{1}|+|Z_{2}|)\Big\}.
$$

By some modification of Lemma A, we obtain the following.\\

\noindent\textbf{Lemma 1}. Let $Q$ be a cycle and $(Z_{1},Z_{2})$ be a nontrivial $(Q,r)$-scheme. Then $|Q|\geq |Z_1|+|Z_2|+|Z_1\cup Z_2|$. Moreover,

$(a1)$ if $r=4$, then $|Q|\geq 2(|Z_1|+|Z_2|)$,

$(a2)$ if $r\geq5$, then $|Q|\geq 2(|Z_1|+|Z_2|)+r-3$,

$(a3)$ if $|Z_2|=1$, then $|Q|\geq2|Z_1|+2r-4$,

$(a4)$ if $|Z_1|=|Z_2|=1$, then $|Q|\geq2r$,

$(a5)$ if $|Z_1|+|Z_2|\geq4$ and $r\geq4$, then $|Q|\geq2(|Z_1|+|Z_2|)+2r-8$.\\

For $(Z_1,Z_2,Z_3)$-schemes we need a result proved in [11, Lemma 1].\\

\noindent\textbf{Lemma B [11]}. Let $Q$ be a cycle and $(Z_{1},Z_{2},Z_3)$ be a nontrivial $(Q,r)$-scheme with $|Z_3|=1$. Then
$$
 |Q|\geq\min\{ 2(|Z_1|+|Z_2|)+3r-10,\frac{1}{2}r(|Z_1|+|Z_2|)\}.
$$

For the special cases $|Z_i|=1$ $(i=2,3)$ and $|Z_i|=1$ $(i=1,2,3)$, we need the following result.\\

\noindent\textbf{Lemma 2}. Let $Q$ be a cycle and $(Z_{1},Z_{2},Z_3)$ be a nontrivial $(Q,r)$-scheme.

(b1) If $r\geq4$, $|Z_1|+|Z_2|\geq 6$ and $|Z_3|=1$, then $|Q|\geq2(|Z_1|+|Z_2|)+3r-12$.

(b2) If $|Z_i|=1$ $(i=2,3)$, then $|Q|\geq2|Z_1|+3r-6$.

(b3) If $|Z_i|=1$ $(i=1,2,3)$, then $|Q|\geq3r$.\\

For $(Z_1,Z_2,Z_3,Z_4)$-schemes we need the next result proved in [11, Lemma 2].\\

\noindent\textbf{Lemma C [11]}. Let $Q$ be a cycle and $(Z_{1},Z_{2},Z_3,Z_4)$ be a nontrivial $(Q,r)$-scheme with $|Z_3|=|Z_4|=1$. Then
$$
|Q|\geq\min\{ 2(|Z_1|+|Z_2|)+4r-14,\frac{1}{2}r(|Z_1|+|Z_2|)\}.
$$

For the special cases $|Z_i|=1$ $(i=2,3,4)$ and $|Z_i|=1$ $(i=1,2,3,4)$, we prove the following.\\

\noindent\textbf{Lemma 3}. Let $Q$ be a cycle and $(Z_{1},Z_{2},Z_3,Z_4)$ be a nontrivial $(Q,r)$-scheme.

(c1) If $r\geq4$, $|Z_1|+|Z_2|\geq8$ and $|Z_3|=|Z_4|=1$, then  $|Q|\geq2(|Z_1|+|Z_2|)+4r-16$.

(c2) If $|Z_i|=1$ $(i=2,3,4)$, then $|Q|\geq2|Z_1|+4r-8$.

(c3) If $|Z_i|=1$ $(i=1,2,3,4)$, then $|Q|\geq4r$.\\

We will show that the main conclusion of Theorem 1 easily follows when $\delta\leq 2\kappa-3$. \\

\noindent\textbf{Lemma 4}. Every 3-connected graph with $\delta\leq 2\kappa-3$ either has a cycle of length at least $4\delta-2\kappa$ or has a dominating cycle.\\

Our constructions of dominating cycles are based on a combination of $Q^{\uparrow}_{1},...,Q^{\uparrow}_{m}$ and $Q^{\downarrow}_{1},...,Q^{\downarrow}_{m}$ (see Definitions A-B). For the case $\delta>3\kappa/3-1$, the existence of $Q^{\downarrow}_{1},...,Q^{\downarrow}_{m}$ follows from the next lemma proved in [11, Lemma 10].\\

\noindent\textbf{Lemma D [11]}. Let $G$ be a 2-connected graph, $A^{\downarrow}$ be an endfragment of $G$ with respect to a minimum cut-set $S$ and let $L$ be a set of independent edges in $\langle S\rangle$. If $\delta > 3\kappa /2 -1$, then $\langle A^{\downarrow}\cup V(L)\rangle$ contains a cycle that uses all the edged in $L$.\\
 
We need also the following result that occurs in [13].\\

\noindent\textbf{Lemma E [13]}. Let $G$ be a hamiltonian graph with $\{v_1,...,v_r\}\subseteq V(G)$ and $d(v_i)\geq r$ $(i=1,...,r)$. Then any two vertices of $G$ are connected by a path of length at least $r$.\\

To prove Theorems 2-5, we need a number of lemmas on total lengths of  $Q^{\uparrow}_{1},...,Q^{\uparrow}_{m}$ and $Q^{\downarrow}_{1},...,Q^{\downarrow}_{m}$ by covering some alternative conditions. In particular, the alternative conditions $|A^{\uparrow}|\leq 3\delta-\kappa-4$ and $|A^{\uparrow}|\geq 3\delta-\kappa-3$ are covered by the following two lemmas.\\

\noindent\textbf{Lemma 5}. Let $G$ be a 3-connected graph and $|A^{\uparrow}|\leq 3\delta-\kappa-4$ for a fragment $A^{\uparrow}$ of $G$. Then $A^{\uparrow}-V^{\uparrow}$ is an independent set.\\

\noindent\textbf{Lemma 6}. Let $G$ be a 4-connected graph and $|A^{\uparrow}|\geq 3\delta-\kappa-3$ for a fragment $A^{\uparrow}$ of $G$. Then either $A^{\uparrow}-V^{\uparrow}$ is an independent set or $|V^{\uparrow}|\geq3\delta-5$.\\

 The next lemma insures $A^{\downarrow}-V^{\downarrow}=\emptyset$ when $|A^{\downarrow}|\leq 2\delta-2\kappa+1$.\\

\noindent\textbf{Lemma 7}. Let $G$ be a 2-connected graph with $|A^{\downarrow}|\leq 2\delta-2\kappa+1$ for an endfragment $A^{\downarrow}$ of $G$. If $\{ Q^{\downarrow}_{1},...,Q^{\downarrow}_{m}\}$ exists, then $A^{\downarrow}-V^{\downarrow}=\emptyset$ and $A^{\downarrow}-V(Q^{\downarrow}_*)=\emptyset$.\\

The alternative conditions $|A^{\downarrow}|\leq 3\delta - 3\kappa+1$ and $|A^{\downarrow}|\geq 3\delta - 3\kappa+2$ are covered by the next two lemmas.\\

\noindent\textbf{Lemma 8}. Let $G$ be a 3-connected graph with $\delta \geq2\kappa -2$ and let $|A^{\downarrow}|\leq 3\delta - 3\kappa+1$ for an endfragment $A^{\downarrow}$ of $G$ with respect to a minimum cut-set S.

$(d1)$ If $f\geq3$, then $A^{\downarrow}-V^{\downarrow}$ is an independent set, where $f=|V^{\downarrow}\cap S|$.

$(d2)$ If $f=2$ and $S\not\subseteq V^{\uparrow}$, then $A^{\downarrow}-V(Q^{\downarrow}_*)$ is an independent set. 

$(d3)$ If $f=2$ and $S\subseteq V^{\uparrow}$, then either $A^{\downarrow}-V^{\downarrow}=\emptyset$ or $|V^{\downarrow}|\geq 2\delta - 2\kappa+3$.\\

\noindent\textbf{Lemma 9}. Let $G$ be a 3-connected graph with $\delta \geq2\kappa  -2$ and let $|A^{\downarrow}|\geq 3\delta - 3\kappa+2$ for an endfragment $A^{\downarrow}$ of $G$ with respect to a minimum cut-set S.

$(e1)$ If $f\geq3$, then either $A^{\downarrow}-V^{\downarrow}$ is an independent set or $|V^{\downarrow}|\geq 3\delta -3\kappa +f$, where $f=|V^{\downarrow}\cap S|$.

$(e2)$ If $f=2$ and $S\not\subseteq V^{\uparrow}$, then either $A^{\downarrow}-V(Q^{\downarrow}_*)$ is an independent set or $|V^{\downarrow}|\geq 3\delta -3\kappa +3$.

$(e3)$ If $f=2$ and $S\subseteq V^{\uparrow}$, then $|V^{\downarrow}|\geq 2\delta -2\kappa +3$.\\ 
 
Some extremal cases are considered in the following lemma.\\

\noindent\textbf{Lemma 10}. Let $G$ be a 2-connected graph and $A^{\uparrow}$ is a fragment of $G$.

$(f1)$ If $A^{\uparrow}-V^{\uparrow}$ is an independent non-empty set, then $|V^{\uparrow}|\geq 2\delta +m-2$.
 
$(f2)$ If $A^{\downarrow}-V^{\downarrow}$ is an independent non-empty set, then $|V^{\downarrow}|\geq 2\delta-2\kappa+2f-m$.

$(f3)$ If $A^{\downarrow}-V(Q^{\downarrow}_*)$ is an independent non-empty set, then $|V^{\downarrow}|\geq 2\delta-2\kappa+5$.\\

For the special case $|A^{\downarrow}|\leq 3\delta - 3\kappa$, we need the following more simple result proved in [11, Lemma 14].\\

\noindent\textbf{Lemma F [11]}. Let $G$ be a 3-connected graph with $\delta\geq2\kappa-2$ and let $|A^{\downarrow}|\leq 3\delta-3\kappa$ for a fragment $A^{\downarrow}$ of $G$. Then $A^{\downarrow}-V^{\downarrow}$ is an independent set.\\

\section{Proofs of lemmas}

\noindent\textbf{Proof of Lemma 1}. Let $x_1,...,x_t$ be the elements of $Z_1\cup Z_2$, occuring on $\overrightarrow{Q}$ in a consecutive order. Consider the segments $I_i=x_i\overrightarrow{Q}x_{i+1}$ $(i=1,...,t)$, where $x_{t+1}=x_1$ and $t=|Z_1\cup Z_2|$. Since $(Z_1,Z_2)$ is a nontrivial $(Q,r)$-scheme, a segment $I_i$ is an $r$-segment if $x_i\in Z_1\cap Z_2$ and is a 2-segment, otherwise. Hence,
$$
|Q|\geq 3|Z_1\cap Z_2|+2(t-|Z_1\cap Z_2|)=|Z_1|+|Z_2|+|Z_1\cup Z_2|.
$$

\textbf{Case 1.} $r=4$.

By Lemma A, $|Q|\geq2(|Z_1|+|Z_2|)$.\\

\textbf{Case 2}. $r\geq5$.

By the hypothesis, $(r-4)(|Z_1|+|Z_2|)\geq2$ which is equivalent to
$$
\frac{1}{2}r(|Z_1|+|Z_2|)\geq2(|Z_1|+|Z_2|)+r-3
$$
and the result follows from Lemma A immediately.\\  

\textbf{Case 3}. $|Z_2|=1$.

Clearly at least two of segments $I_1,...,I_t$ are $r$-segments. If $Z_2\subseteq Z_1$, then $t=|Z_1|$ and $|C|\geq2r+2(t-2)=2|Z_1|+2r-4$. If $Z_2\not\subseteq Z_1$, then $t=|Z_1|+1$ and $|Q|\geq2r+2(t-1)>2|Z_1|+2r-4$.\\

\textbf{Case 4}. $|Z_1|=|Z_2|=1$.

The result is trivial.\\

\textbf{Case 5}. $r\geq4$ and $|Z_1|+|Z_2|\geq4$.

By the hypothesis, $(r-4)(|Z_1|+|Z_2|-4)\geq0$ which is equivalent to
$$
\frac{1}{2}r(|Z_1|+|Z_2|)\geq2(|Z_1|+|Z_2|)+2r-8
$$
and the result follows from Lemma A, immediately. \qquad  $\Delta$\\

\noindent\textbf{Proof of Lemma 2}. Let $x_1,...,x_t$ be the elements of $Z_1\cup Z_2\cup Z_3$, occuring on $\overrightarrow{Q}$ in a consecutive order. Consider the segments $I_i=x_i\overrightarrow{Q}x_{i+1}$ $(i=1,...,t)$, where $x_{t+1}=x_1$ and $t=|Z_1\cup Z_2\cup Z_3|$. \\  

\textbf{Case 1}. $r\geq4$, $|Z_1|+|Z_2|\geq6$ and $|Z_3|=1$.

By the hypothesis, $(r-4)(|Z_1|+|Z_2|-6)\geq0$ which is equivalent to
$$
\frac{1}{2}r(|Z_1|+|Z_2|)\geq2(|Z_1|+|Z_2|)+3r-12
$$
and the result follows from Lemma B, immediately. \\

\textbf{Case 2}. $|Z_2|=|Z_3|=1$.

Clearly at least three of segments $I_1,...,I_t$ are $r$-segments. If $Z_2\cup Z_3\subseteq Z_1$, then $t=|Z_1|$ and $|Q|\geq3r+2(t-3)\geq2|Z_1|+3r-6$. If $Z_2\subseteq Z_1$ and $Z_3\not\subseteq Z_1$ or $Z_3\subseteq Z_1$ and $Z_2\not\subseteq Z_1$, then $t=|Z_1|+1$ and
$$
|Q|\geq3r+2(t-3)=2|Z_1|+3r-4>2|Z_1|+3r-6.
$$

Finally, if $Z_2\not\subseteq Z_1$ and $Z_3\not\subseteq Z_1$, then $t=|Z_1|+2$ and 
$$
|Q|\geq3r+2(t-3)=2|Z_1|+3r-2>2|Z_1|+3r-6.
$$

\textbf{Case 3}. $|Z_1|=|Z_2|=|Z_3|=1$.

The result is trivial.   \qquad  $\Delta$\\

\noindent\textbf{Proof of Lemma 3}. Let $x_1,...,x_t$ be the elements of $Z_1\cup Z_2\cup Z_3\cup Z_4$, occuring on $\overrightarrow{Q}$ in a consecutive order. Consider the segments $I_i=x_i\overrightarrow{Q}x_{i+1}$ $(i=1,...,t)$, where $x_{t+1}=x_1$ and $t=|Z_1\cup Z_2\cup Z_3\cup Z_4|$. \\ 

\textbf{Case 1}. $r\geq4$, $|Z_1|+|Z_2|\geq8$ and $|Z_3|=|Z_4|=1$.

By the hypothesis, $(r-4)(|Z_1|+|Z_2|-8)\geq0$ which is equivalent to
$$
\frac{1}{2}r(|Z_1|+|Z_2|)\geq2(|Z_1|+|Z_2|)+4r-16
$$
and the result follows from Lemma C, immediately. \\

\textbf{Case 2}. $|Z_i|=1$ (i=2,3,4).

We can argue as in proof of Lemma 2 (Case 2).\\

\textbf{Case 3}. $|Z_i|=1$ (i=1,2,3,4).

The result is trivial.   \qquad  $\Delta$\\

\noindent\textbf{Proof of Lemma 4}. If $G$ has a dominating cycle, then we are done. Otherwise, using Theorem M and the fact that $\delta\leq2\kappa-3$, we obtain $c\geq 3\delta-3\geq4\delta-2\kappa$. \quad  $\Delta$\\ 

\noindent\textbf{Proofs of Lemmas 5-6}. Let $V^{\uparrow}$ be as defined in Definition A and let $P$ be a longest path in $\langle A^{\uparrow}-V^{\uparrow}\rangle$. If $|V(P)|\leq1$, then clearly $A^{\uparrow}-V^{\uparrow}$ is an independent set. Let $|V(P)|\geq2$. To prove Lemma 5, we can argue exactly as in [11, proof of Lemma 12]. To prove Lemma 6, we can argue exactly as in [11, proof of Lemma 13].\qquad   $\Delta$\\

\noindent\textbf{Proof of Lemma 7}. Let $A^{\downarrow}$ is defined with respect to a minimum cut-set $S$ and let $Q^{\downarrow}_{1},...,Q^{\downarrow}_{m}$, $V^{\downarrow}$  and $C$ be as defined in Definition B. Put $f=|V^{\downarrow}\cap S|$. Assume w.l.o.g. that $F(Q^{\uparrow}_i)=v_i$ and $L(Q^{\uparrow}_i)=w_i$ $(i=1,...,m)$. Form a cycle $C^{\downarrow}$ by deleting $Q^{\uparrow}_{1},...,Q^{\uparrow}_{m}$ from $C$ and adding extra edges $v_1w_1,v_2w_2,...,v_mw_m$, respectively. Clearly $V^{\downarrow}=V(C^{\downarrow})$.  Conversely, suppose that $A^{\downarrow}-V^{\downarrow}\neq\emptyset$, and choose a connected component $H$ of $\langle A^{\downarrow}-V^{\downarrow}\rangle$. Put $M=N(V(H))$. If $|M|=\kappa$, then $M$ is a minimum cut-set of $G$ with $M\subset A^{\downarrow}\cup S$ and $M\neq S$, contradicting the definition of $A^{\downarrow}$. Let $|M|\geq\kappa+1$. Clearly $|M\cap (S-V^{\downarrow})|\leq \kappa-f$. Since  $M\subseteq V(C^{\downarrow})\cup S$, we have 
$$
|M\cap V(C^{\downarrow})|=|M|-|M\cap (S-V^{\downarrow})|\geq f+1.
$$

Recalling that $|V(C^{\downarrow})\cap S|=f$, we have $x^+\notin S$ for some $x\in M\cap V(C^{\downarrow})$. Then by deleting $xx^+$ from $C^{\downarrow}$ and adding an edge $xy$ with $y\in V(H)$, we get a path $R=x^+\overrightarrow{C}^{\downarrow}xy$ with endvertices in $A^{\downarrow}$. Choose a maximal path $R^*=\xi R^*\eta$ in $\langle A^{\downarrow}\cup V^{\downarrow}\rangle$ containing $R$ as a subpath with $\xi ,\eta\in A^{\downarrow}$ and put
$$
d_1=|N(\xi)\cap V(R^*)|, \quad  d_2=|N(\eta)\cap V(R^*)|.
$$

Since $R^{*}$ is extreme, we have $d_i\geq\delta-\kappa+f$ $(i=1,2)$ and
$$
|V(R^*)|\leq|A^{\downarrow}|+f\leq2\delta-2\kappa+f+1\leq d_1+d_2-f+1,
$$

\noindent implying that $d_1+d_2\geq|V(R^*)|+f-1\geq|V(R^*)|+m$. By standard arguments, there is a vertex $\lambda\in V(\overrightarrow{R}^*)$ such that $\xi\lambda^+\in E(G)$, $\eta\lambda\in E(G)$ and $\lambda\lambda^+\not\in\{v_1w_1,...,v_mw_m\}$. Then by deleting $\lambda\lambda^+$ and adding $\xi\lambda^+$, $\eta\lambda$, we can form a new collection of paths instead of $Q^{\downarrow}_{1},...,Q^{\downarrow}_{m}$, contrary to $(*2)$. So, $A^{\downarrow}-V^{\downarrow}=\emptyset$. By the same way it can be shown that $A^{\downarrow}-V(Q^{\downarrow}_*)=\emptyset$.  \qquad   $\Delta$ \\

\noindent\textbf{Proof of Lemma 8}. Let $Q^{\downarrow}_{1},...,Q^{\downarrow}_{m}$ and $V^{\downarrow}_{1},...,V^{\downarrow}_{m},V^{\downarrow}$ are as defined in Definition B. The existence of $Q^{\downarrow}_{1},...,Q^{\downarrow}_{m}$ follows from Lemma D. Clearly $f\geq2m$. If $|A^{\downarrow}|\leq2\delta-2\kappa+1$, then by Lemma 7, we are done. Let $|A^{\downarrow}|\geq2\delta-2\kappa+2$\\

\textbf{Case 1}. $f\geq3$.

Suppose first that $\delta-\kappa\leq1$. Combining this with $\delta\geq2\kappa-2$, we get $f=\kappa=3$, $\delta=4$, $m=1$ and $|A^{\downarrow}|\leq4$. If $A^{\downarrow}-V^{\downarrow}$ is an independent set, then we are done. Otherwise, $|V^{\downarrow}|\leq |A^{\downarrow}|+f-2\leq5$. Choose a vertex $z\in A^{\downarrow}-V^{\downarrow}$. If $z$ and $Q^{\downarrow}_{1}$ are connected by at most three paths having no vertex other than $z$ in common, then $G$ has a cut-set $S^{\prime}$ of order three with $S^{\prime}\subset A^{\downarrow}\cup S$ and $S^{\prime}\neq S$, contradicting the definition of $A^{\downarrow}$. Otherwise, it is easy to see that $|V^{\downarrow}|\geq7$, a contradiction. So, we can assume that $\delta-\kappa\geq2$. Let $P=y_1...y_p$ be a longest path in $\langle A^{\downarrow}-V^{\downarrow}\rangle$. By the hypothesis,
$$
p+|V^{\downarrow}|-f\leq|A^{\downarrow}|\leq3\delta-3\kappa+1,
$$
implying that
$$
|V^{\downarrow}|\leq3\delta-3\kappa+f-p+1. \eqno{(1)}
$$

Put
$$
Z_1=N(y_1)\cap V^{\downarrow},\quad Z_2=N(y_p)\cap V^{\downarrow},
$$
$$
Z_{1,i}=Z_1\cap V^{\downarrow}_{i},\quad Z_{2,i}=Z_2\cap V^{\downarrow}_{i}\quad (i=1,...,m).
$$

Clearly $Z_1=\cup^{m}_{i=1}Z_{1,i}$ and $Z_2=\cup^{m}_{i=1}Z_{2,i}$. If $p\leq1$, then $A^{\downarrow}-V^{\downarrow}$ is independent and we are done. Let $p\geq2$.\\

\textbf{Case 1.1}. $p=2$.

Clearly $|Z_i|\geq\delta-\kappa+f-1$ $(i=1,2)$. For each $i\in \{1,...,m\}$, form a cycle $C^{\downarrow}_i$ by adding to $Q^{\downarrow}_i$ an extra path of length 3 connecting $F(Q^{\downarrow}_i)$ and $L(Q^{\downarrow}_i)$. Since $\{Q^{\downarrow}_{1},...,Q^{\downarrow}_{m}\}$ is $(*2)$-extreme, $(Z_{1,i},Z_{2,i})$ is a $(C^{\downarrow}_{i},3)$-scheme. If $(Z_{1,i},Z_{2,i})$ is a nontrivial $(C^{\downarrow}_{i},3)$-scheme, then by Lemma 1, $|C^{\downarrow}_{i}|\geq2|Z_{1,i}|+|Z_{2,i}|$. Otherwise, it can be checked easily. Hence,
$$
|V^{\downarrow}|=\sum^{m}_{i=1}|V^{\downarrow}_{i}|=\sum^{m}_{i=1}(|C^{\downarrow}_{i}|-2)\geq\sum^{m}_{i=1}|C^{\downarrow}_{i}|-2m
$$
$$
\geq2|Z_1|+|Z_2|-2m\geq(3\delta-3\kappa+f)+2f-2m-3.
$$

Since $f\geq2m$ and $f\geq3$, we have $2f-2m-3\geq0$, implying that $|V^{\downarrow}|\geq3\delta-3\kappa+f$, contrary to (1).\\

\textbf{Case 1.2.} $p=3$.

Clearly $|Z_i|\geq\delta-\kappa+f-2$ $(i=1,2)$. For each $i\in \{1,...,m\}$ form a cycle $C^{\downarrow}_i$ by adding to $Q^{\downarrow}_i$ an extra path of length 4, connecting $F(Q^{\downarrow}_i)$ and $L(Q^{\downarrow}_i)$. Since $\{Q^{\downarrow}_{1},...,Q^{\downarrow}_{m}\}$ is $(*2)$-extreme, $(Z_{1,i},Z_{2,i})$ is a $(C^{\downarrow}_{i},4)$-scheme. If $(Z_{1,i},Z_{2,i})$ is a nontrivial $(C^{\downarrow}_{i},4)$-scheme, then by $(a1)$ (see Lemma 1), $|C^{\downarrow}_{i}|\geq2(|Z_{1,i}|+|Z_{2,i}|)$. Otherwise, it can be checked easily. Hence,
$$
|V^{\downarrow}|=\sum^{m}_{i=1}|V^{\downarrow}_{i}|=\sum^{m}_{i=1}(|C^{\downarrow}_{i}|-3)\geq\sum^{m}_{i=1}|C^{\downarrow}_{i}|-3m
$$
$$
\geq2(|Z_1|+|Z_2|)-3m\geq(3\delta-3\kappa+f)+3f-3m-6.
$$

If $f=3$, then $m=1$ and $3f-3m-6=0$. If $f=4$, then $m\leq2$ and $3f-3m-6\geq0$. Finally, if $f\geq5$, then again $3f-3m-6\geq f+m-6\geq0$. So, in any case, $|V^{\downarrow}|\geq3\delta-3\kappa+f$, contrary to (1).\\

\textbf{Case 1.3.} $p\geq4$.

Let $w_1,w_2,...,w_s$ be the elements of $(N(y_p)\cap V(P))^+$ occuring on $\overrightarrow{P}$ in a consecutive order, where $w_s=y_p$. Put $P_0=w^-_1\overrightarrow{P}w_s$ and $p_0=|V(P_0)|$. For each $w_i$ $(i\in \{1,...,s\})$ there is a path $y_1\overrightarrow{P}w^{-}_{i}w_s\overleftarrow{P}w_i$ in $\langle V(P)\rangle$ of length $p$ connecting $y_1$ and $w_i$. Hence, we can assume w.l.o.g. that for each $i\in \{1,...,s\}$, $P$ is chosen such that 
$$
|Z_1|\geq|N(w_i)\cap V^{\downarrow}|,\quad N(w_i)\cap V(P)\subseteq V(P_0).
$$ 

Choose $w\in\lbrace w_1,...,w_s\rbrace$ as to maximize $|N(w_i)\cap V^{\downarrow}|$, $i=1,...,s$. Set
$$
Z_3=N(w)\cap V^{\downarrow}, \quad Z_{3,i}=Z_3\cap V^{\downarrow}_{i} \quad (i=1,...,m).
$$

Clearly 
$$
|Z_1|\geq|Z_3|\geq|Z_2|\geq\delta-\kappa+f-p_0+1.    \eqno{(2)}   
$$
 
\textbf{Claim 1.} $|V^{\downarrow}_{i}|\geq2(|Z_{1,i}|+|Z_{3,i}|)-2$ $(i=1,...,m)$. 

Proof. For each $i\in \{1,...,m\}$, form a cycle $C^{\downarrow}_{i}$ by adding to $Q^{\downarrow}_{i}$ an extra path of length $p+1$ connecting $F(Q^{\downarrow}_i)$ and $L(Q^{\downarrow}_i)$. Since $\{Q^{\downarrow}_{1},...,Q^{\downarrow}_{m}\}$ is $(*2)$-extreme, $(Z_{1,i},Z_{3,i})$ is a $(C^{\downarrow}_{i},p+1)$-scheme. If $(Z_{1,i},Z_{3,i})$ is a nontrivial $(C^{\downarrow}_{i},p+1)$-scheme, then by $(a2)$ (see Lemma 1), $|C^{\downarrow}_{i}|\geq2(|Z_{1,i}|+|Z_{3,i}|)+p-2$. Otherwise, it can be checked easily. Observing that $|V^{\downarrow}_{i}|=|C^{\downarrow}_{i}|-p$ $(i=1,...,m)$, we get the desired result immediately.  \qquad $\Delta$\\

By Claim 1,
$$
|V^{\downarrow}|=\sum^{m}_{i=1}|V^{\downarrow}_{i}|\geq2(|Z_1|+|Z_3|)-2m\geq4(\delta-\kappa+f-p_0+1)-2m
$$
$$
=(3\delta-3\kappa+f)+(\delta-\kappa-2)+4(m-p_0+1)+3(f-2m)+2
$$
$$
>(3\delta-3\kappa+f)+4(m-p_0+1).
$$

If $p_0\leq m+1$, then clearly $|V^{\downarrow}|\geq3\delta-3\kappa+f$, contrary to (1). So, we can assume that
$$
p_0 \geq m+2.   \eqno{(3)}
$$

Further, assume that $\delta-\kappa+f-|Z_3|\leq2$, implying that $|Z_1|\geq|Z_3|\geq\delta-\kappa+f-2$. By Claim 1,  
$$
|V^{\downarrow}|=\sum^{m}_{i=1}|V^{\downarrow}_{i}|\geq2(|Z_1|+|Z_3|)-2m\geq4(\delta-\kappa+f-2)-2m
$$
$$
=(3\delta-3\kappa+f)+(\delta-\kappa-2)+3f-2m-6\geq3\delta-3\kappa+f,
$$
contrary to (1). Now let $\delta-\kappa+f-|Z_3|\geq3$. By the choice of $w$, 
$$
|N(w_i)\cap V(P_0)|\geq\delta-\kappa+f-|Z_3|\geq3 \quad (i=1,...,s).
$$ 

In particular, for $i=s$, we have $s\geq\delta-\kappa+f-|Z_3|$. By Lemma E, in $\langle V(P_0)\rangle$ any two vertices are joined by a path of length at least $\delta-\kappa+f-|Z_3|$. Observing that $p\geq s+1\geq \delta-\kappa+f-|Z_3|+1$ and combining it with (1), we get 
$$
|V^{\downarrow}|\leq2\delta-2\kappa+|Z_3|.   \eqno{(4)}
$$

For each $i\in \{1,...,m\}$, form a cycle $C^{\downarrow}_{i}$ by adding to $Q^{\downarrow}_{i}$ an extra path of length $r$, where $r=\delta-\kappa+f-|Z_3|+2\geq5$, connecting $F(Q^{\downarrow}_i)$ and $L(Q^{\downarrow}_i)$. Let $G^{\prime}$ be the subgraph of $G$ obtained from $G$ by deleting all vertices of $S-V^{\downarrow}$. Clearly $G^{\prime}$ is $f$-connected. If $m=1$, then $f\geq3=m+2$. If $m\geq2$, then again $f\geq2m\geq m+2$. So, in any case, $G^{\prime}$ is $(m+2)$-connected. Since $|V^{\downarrow}|>f\geq m+2$ and $|V(P_0)|=p_0\geq m+2$ (by (3)), there are vertex disjoint paths $R_1,...,R_{m+2}$ in $G^{\prime}$ connecting $V(P_0)$ and $V^{\downarrow}$. Let $F(R_i)\in V(P_0)$ and $L(R_i)\in V^{\downarrow}$ $(i=1,...,m+2)$. Then we can assume w.l.o.g. that either $L(R_i)\in V^{\downarrow}_{1}$ $(i=1,2)$ and $L(R_i)\in V^{\downarrow}_{2}$ $(i=3,4)$ or $L(R_i)\in V^{\downarrow}_{1}$ $(i=1,2,3)$.  \\

\textbf{Claim 2.} Let $\tau\in\{1,...,m\}$ and $a,b\in\{1,...,m+2\}$. If $L(R_i)\in V^{\downarrow}_{\tau}$ $(i=a,b)$, then $|V^{\downarrow}_{\tau}|\geq2(|Z_{1,\tau}|+|Z_{3,\tau}|)+r-7$.

Proof. Suppose first that $|Z_{1,\tau}|\leq2$, $|Z_{3,\tau}|\leq2$. Clearly $(\{L(R_a)\},\{L(R_b)\})$ is a nontrivial $(C^{\downarrow}_{\tau},r)$-scheme. By $(a4)$ (see Lemma 1), $|C^{\downarrow}_{\tau}|\geq2r$, implying that
$$
|V^{\downarrow}_{\tau}|=|C^{\downarrow}_{\tau}|-r+1\geq r+1\geq 2(|Z_{1,\tau}|+|Z_{3,\tau}|)+r-7.
$$

Now let $|Z_{1,\tau}|\geq3$ and $|Z_{3,\tau}|\leq2$. If either $R_a$ or $R_b$, say $R_a$, has no common vertices with $y_1\overrightarrow{P}w^-_1$, then $(Z_{1,\tau},\{L(R_a)\})$ is a nontrivial $(C^{\downarrow}_{\tau},r)$-scheme. Otherwise, let $t$ be the smallest integer such that $y_t\in V(R_a\cup R_b)$. Assume w.l.o.g. that $y_t\in V(R_b)$. Then due to $y_1\overrightarrow{P}y_tR_bF(R_b)$, we again can state that $(Z_{1,\tau},\{L(R_a)\})$ is a nontrivial $(C^{\downarrow}_{\tau},r)$-scheme. By $(a3)$ (see Lemma 1), $|C^{\downarrow}_{\tau}|\geq2|Z_{1,\tau}|+2r-4$, implying that
$$
|V^{\downarrow}_{\tau}|=|C^{\downarrow}_{\tau}|-r+1\geq 2|Z_{1,\tau}|+r-3\geq 2(|Z_{1,\tau}|+|Z_{3,\tau}|)+r-7.
$$

Finally, we suppose that $|Z_{1,\tau}|\geq3$ and $|Z_{3,\tau}|\geq3$. It means that $(Z_{1,\tau},Z_{3,\tau})$ is a nontrivial $(C^{\downarrow}_{\tau},r)$-scheme. By $(a5)$ (see Lemma 1), $|C^{\downarrow}_{\tau}|\geq2(|Z_{1,\tau}|+|Z_{3,\tau}|)+2r-8$, implying that
$$
|V^{\downarrow}_{\tau}|=|C^{\downarrow}_{\tau}|-r+1\geq 2(|Z_{1,\tau}|+|Z_{3,\tau}|)+r-7. \qquad \Delta\\
$$

\textbf{Claim 3.} Let $\tau\in\{1,...,m\}$ and $a,b,c\in\{1,...,m+2\}$. If $L(R_i)\in V^{\downarrow}_{\tau}$ $(i=a,b,c)$, then $|V^{\downarrow}_{\tau}|\geq2(|Z_{1,\tau}|+|Z_{3,\tau}|)+2r-11$.

Proof. Suppose first that $|Z_{1,\tau}|\leq3$ and $|Z_{3,\tau}|\leq3$. Clearly $\{L(R_a)\}$, $\{L(R_b)\}$ and $\{L(R_c)\}$ form a nontrivial $(C^{\downarrow}_{\tau},r)$-scheme. By $(b3)$ (see Lemma 2), $|C^{\downarrow}_{\tau}|\geq3r$ and therefore,
$$
|V^{\downarrow}_{\tau}|=|C^{\downarrow}_{\tau}|-r+1\geq 2r+1\geq 2(|Z_{1,\tau}|+|Z_{3,\tau}|)+2r-11.
$$

Now let $|Z_{1,\tau}|\geq4$ and $|Z_{3,\tau}|\leq3$. If at least two of $R_a$, $R_b$, $R_c$, say $R_a$ and $R_b$, have no common vertices with $y_1\overrightarrow{P}w^-_1$, then $Z_{1,\tau}$, $\{L(R_a)\}$ and $\{L(R_b)\}$ form a nontrivial $(C^{\downarrow}_{\tau},r)$-scheme. Otherwise, let $t$ be the smallest integer such that $y_t\in V(R_a\cup R_b\cup R_c)$. Assume w.l.o.g. that $y_t\in V(R_c)$. Then due to $y_1\overrightarrow{P}y_tR_cF(R_c)$, we again can state that $Z_{1,\tau}$, $\{L(R_a)\}$ and $\{L(R_b)\}$ form a nontrivial $(C^{\downarrow}_{\tau},r)$-scheme. By $(b2)$ (see Lemma 2), $|C^{\downarrow}_{\tau}|\geq2|Z_{1,\tau}|+3r-6$, implying that
$$
|V^{\downarrow}_{\tau}|=|C^{\downarrow}_{\tau}|-r+1\geq 2|Z_{1,\tau}|+2r-5\geq 2(|Z_{1,\tau}|+|Z_{3,\tau}|)+2r-11.
$$

Finally, let $|Z_{1,\tau}|\geq4$ and $|Z_{3,\tau}|\geq4$. If at least two of $R_a$, $R_b$, $R_c$, say $R_a$ and $R_b$, have no common vertices with $y_1\overrightarrow{P}w^-_1$, then   choosing one of $R_a$, $R_b$, say $R_a$, satisfying $F(R_a)\neq w$, we conclude that $Z_{1,\tau}$, $Z_{3,\tau}$ and $\{L(R_a)\}$ form a nontrivial $(C^{\downarrow}_{\tau},r)$-scheme. Otherwise, let $t$ be the smallest integer such that $y_t\in V(R_a\cup R_b\cup R_c)$. Assume w.l.o.g. that $y_t\in V(R_c)$. Choosing one of $R_a$, $R_b$, say $R_a$, satisfying $F(R_a)\neq w$, we again can state that $(Z_{1,\tau},Z_{3,\tau},\{L(R_a)\})$ is a nontrivial $(C^{\downarrow}_{\tau},r)$-scheme. By $(b1)$ (see Lemma 2), $|C^{\downarrow}_{\tau}|\geq2(|Z_{1,\tau}|+|Z_{3,\tau}|)+3r-12$, implying that
$$
|V^{\downarrow}_{\tau}|=|C^{\downarrow}_{\tau}|-r+1\geq 2(|Z_{1,\tau}|+|Z_{3,\tau}|)+2r-11. \qquad \Delta
$$

\textbf{Claim 4.} Let $\tau\in\{1,...,m\}$ and $a,b,c,d\in\{1,...,m+3\}$. If $L(R_i)\in V^{\downarrow}_{\tau}$ $(i=a,b,c,d)$, then $|V^{\downarrow}_{\tau}|\geq2(|Z_{1,\tau}|+|Z_{3,\tau}|)+3r-15$.

Proof. Suppose, that $|Z_{1,\tau}|\leq4$, $|Z_{3,\tau}|\leq4$. Clearly $\{L(R_a)\}$, $\{L(R_b)\}$, $\{L(R_c)\}$ and $\{L(R_d)\}$ form a nontrivial $(C^{\downarrow}_{\tau},r)$-scheme. By $(c3)$ (see Lemma 3), $|C^{\downarrow}_{\tau}|\geq4r$, implying that
$$
|V^{\downarrow}_{\tau}|=|C^{\downarrow}_{\tau}|-r+1\geq 3r+1\geq 2(|Z_{1,\tau}|+|Z_{3,\tau}|)+3r-15.
$$

Now let $|Z_{1,\tau}|\geq5$ and $|Z_{3,\tau}|\leq4$. If at least three of $R_a$, $R_b$, $R_c$, $R_d$ say $R_a$, $R_b$, $R_c$, do not intersect $P_0$, then
$Z_{1,\tau}$, $\{L(R_a)\}$, $\{L(R_b)\}$ and $\{L(R_c)\}$ form a nontrivial $(C^{\downarrow}_{\tau},r)$-scheme. Otherwise, let $t$ be the smallest integer such that $y_t\in V(R_a\cup R_b\cup R_c\cup R_d)$. Assume w.l.o.g. that $y_t\in V(R_d)$. Then due to $y_1\overrightarrow{P}y_tR_dF(R_d)$, we again can state that $Z_{1\tau}$, $\{L(R_a)\}$, $\{L(R_b)\}$ and $\{L(R_c)\}$ form a nontrivial $(C^{\downarrow}_{\tau},r)$-scheme. By $(c2)$ (see Lemma 3), $|C^{\downarrow}_{\tau}|\geq2|Z_{1,\tau}|+4r-8$, implying that
$$
|V^{\downarrow}_{\tau}|=|C^{\downarrow}_{\tau}|-r+1\geq 2|Z_{1,\tau}|+3r-7\geq 2(|Z_{1,\tau}|+|Z_{3,\tau}|)+3r-15.
$$

Finally, let $|Z_{1,\tau}|\geq5$ and $|Z_{3,\tau}|\geq5$. If at least three of $R_a$, $R_b$, $R_c$, $R_d$ say $R_a$, $R_b$, $R_c$, do not intersect $P_0$, then   choosing two of $R_a$, $R_b$, $R_c$, say $R_a$, $R_b$, such that $F(R_a)\neq w$ and $F(R_b)\neq w$, we conclude that $Z_{1,\tau}$, $Z_{3,\tau}$, $\{L(R_a)\}$ and $\{L(R_b)\}$ form a nontrivial $(C^{\downarrow}_{\tau},r)$-scheme. Otherwise, let $t$ be the smallest integer such that $y_t\in V(R_a\cup R_b\cup R_c\cup R_d)$. Assume w.l.o.g. that $y_t\in V(R_d)$. Then choosing two of $R_a$, $R_b$, $R_c$, say $R_a$, $R_b$, satisfying $F(R_a)\neq w$ and $F(R_a)\neq w$, we again can state that $Z_{1\tau}$, $Z_{3\tau}$, $\{L(R_a)\}$ and $\{L(R_b)\}$ form a nontrivial $(C^{\downarrow}_{\tau},r)$-scheme. By $(c1)$ (see Lemma 3), $|C^{\downarrow}_{\tau}|\geq2(|Z_{1,\tau}|+|Z_{3,\tau}|)+4r-16$, implying that
$$
|V^{\downarrow}_{\tau}|=|C^{\downarrow}_{\tau}|-r+1\geq 2(|Z_{1,\tau}|+|Z_{3,\tau}|)+3r-15. \qquad \Delta
$$

\textbf{Case 1.3.1}. $L(R_i)\in V^{\downarrow}_{1}$ $(i=1,2)$ and $L(R_i)\in V^{\downarrow}_{2}$ $(i=3,4)$. 

It follows that $f\geq4$. By Claim 2, $|V^{\downarrow}_{i}|\geq2(|Z_{1,i}|+|Z_{3,i}|)+r-7$ $(i=1,2)$. By Claim 1, $|V^{\downarrow}_{i}|\geq2(|Z_{1,i}|+|Z_{3,i}|)-2$ $(i=3,...,m)$. By summing, we get           \\
$$
|V^{\downarrow}|=|V^{\downarrow}_{1}|+|V^{\downarrow}_{2}|+\sum^{m}_{i=3}|V^{\downarrow}_{i}|\geq2(|Z_1|+|Z_3|)+2r-14-2(m-2)
$$
$$
=(2\delta-2\kappa+|Z_3|+1)+f+|Z_3|-7\geq2\delta-2\kappa+|Z_3|+1,
$$
contrary to (4).\\

\textbf{Case 1.3.2}. $L(R_i)\in V^{\downarrow}_{1}$ $(i=1,2,3)$. 

By Claim 3, $|V^{\downarrow}_{1}|\geq2(|Z_{1,1}|+|Z_{3,1}|)+2r-11$. By Claim 1, $|V^{\downarrow}_{i}|\geq2(|Z_{1,i}|+|Z_{3,i}|)-2$ $(i=2,...,m)$. By summing, we get
$$
|V^{\downarrow}|\geq2(Z_1|+|Z_3|)+2r-11-2(m-1)
$$
$$
\geq(2\delta-2\kappa+|Z_3|+1)+f+|Z_3|-6\geq2\delta-2\kappa+|Z_3|+1,
$$
contrary to (4).\\

\textbf{Case 2}. $f=2$ and $S\not\subseteq V^{\uparrow}$.

Choose a vertex $z\in S-V^{\uparrow}$. Put $F(Q^{\downarrow}_1)=v_1$ and $L(Q^{\downarrow}_1)=w_1$. If there is a cut-vertex $x$ in $\langle A^{\downarrow}\cup\{v_1,w_1,z\}\rangle$ that separates $z$ and $Q^{\downarrow}_{1}$, then $S^{\prime }=\{S-z\}\cup \{x\}$  is a minimum cut-set of $G$ other than $S$ with $S^{\prime}\subset A^{\downarrow}\cup S$, contradicting the definition of $A^{\downarrow}$. Otherwise, the existence of $Q^{\downarrow}_{*}$ (see Definition B) follows easily. Put $g=|V(Q^{\downarrow}_{*})\cap S|$. By the definition, $g\geq3$. As in Case 1, we can assume that $\delta-\kappa\geq2$. Let $P=y_1...y_p;Z_1;Z_2$ be as defined in Case 1. If $p\leq1$, then $A^{\downarrow}-V(Q^{\downarrow}_{*})$ is independent. Let $p\geq2$. Put $S_0=(V(Q^{\downarrow}_{*})\cap S)-\{v_1,w_1\}$.\\

\textbf{Case 2.1.} $p=2$.

Clearly $|Z_i|\geq\delta-\kappa+g-1\geq4$ $(i=1,2)$. Form a cycle $C^{\downarrow}_{1}$ by adding to $Q^{\downarrow}_{*}$ an extra path of length 3 with endvertices $v_1,w_1$. Let $I_1,...,I_t$ be the segments of $C^{\downarrow}_{1}$ having only their ends in common with $Z_1\cup Z_2$ and having at least one inner vertex in common with $S$. \\

\textbf{Case 2.1.1.} $S_0\not\subseteq V(I_i)-\{F(I_i),L(I_i)\}$ $(i=1,...,t)$.

It is easy to see that $(Z_1,Z_2)$ is a nontrivial $(C^{\downarrow}_{1},3)$-scheme. By Lemma 1,
$$
|V^{\downarrow}|\geq|V(Q^{\downarrow}_{*})|=|C^{\downarrow}_{1}|-2\geq2|Z_1|+|Z_2|-2\geq3\delta-3\kappa+3g-5.
$$

Hence, 
$$
|A^{\downarrow}|\geq|V^{\downarrow}|-g+2\geq(3\delta-3\kappa+2)+2g-5>3\delta-3\kappa+2,
$$
contradicting the hypothesis.\\

\textbf{Case 2.1.2.} $S_0\subseteq V(I_i)-\{F(I_i),L(I_i)\}$ for some $i\in \{1,...,t\}$, say $i=1$.

Let $I_1=x\overrightarrow{I_1}y$. Denote by $R$ a longest path connecting $x$ and $y$ and passing through $V(P)$. Clearly $|R|=3$ if $x$ and $y$ belong to different $Z_1,Z_2$ and $|R|=2$, otherwise. Form a new cycle $C^{\downarrow}_{2}$ by deleting $I_1$ from $C^{\downarrow}_{1}$ and adding $R$. Clearly $(Z_1,Z_2)$ is a nontrivial $(C^{\downarrow}_{2},3)$-scheme and we can reach a contradiction as in Case 2.1.1.\\

\textbf{Case 2.2.} $p=3$.

Clearly $|Z_i|\geq\delta-\kappa+g-2\geq3$ $(i=1,2)$. Form a cycle $C^{\downarrow}_{1}$ by adding to $Q^{\downarrow}_{*}$ an extra path of length 4 with endvertices $v_1,w_1$. Let $I_1,...,I_t$ be as defined in Case 2.1.\\

\textbf{Case 2.2.1.} $S_0\not\subseteq V(I_i)-\{F(I_i),L(I_i)\}$ $(i=1,...,t)$.

It is easy to see that $(Z_1,Z_2)$ is a nontrivial $(C^{\downarrow}_{1},4)$-scheme. By Lemma 1,
$$
|V^{\downarrow}|\geq|V(Q^{\downarrow}_{*})|=|C^{\downarrow}_{1}|-3\geq2|Z_1|+|Z_2|-3\geq4\delta-4\kappa+4g-11.
$$

Hence,
$$
|A^{\downarrow}|\geq|V^{\downarrow})|-g+3\geq(3\delta-3\kappa+2)+\delta-\kappa+3g-10>3\delta-3\kappa+2,
$$
contradicting the hypothesis.\\

\textbf{Case 2.2.2.} $S_0\subseteq V(I_i)-\{F(I_i),L(I_i)\}$ for some $i\in \{1,...,t\}$, say $i=1$.

Let $I_1=x\overrightarrow{I_1}y$. Denote by $R$ a longest path connecting $x$ and $y$ and passing through $V(P)$. Clearly $|R|=4$ if $x$ and $y$ belong to different $Z_1,Z_2$ and $|R|\geq2$, otherwise. Form a new cycle $C^{\downarrow}_{2}$ by deleting $I_1$ from $C^{\downarrow}_{1}$ and adding $R$. Clearly $(Z_1,Z_2)$ is a nontrivial $(C^{\downarrow}_{2},4)$-scheme and we can reach a contradiction as in Case 2.2.1.\\

\textbf{Case 2.3.} $p\geq4$.

Let $P_0,p_0,w$ and $Z_3$ are as defined in Case 1.3. By the definition, $|Z_1|\geq|Z_3|\geq|Z_2|$ and $p_0\geq2$. If $p_0=2$, then $|Z_1|\geq|Z_2|\geq\delta-\kappa+g-1$ and we can argue as in case $p=2$. If $p_0=3$, then $|Z_1|\geq|Z_2|\geq\delta-\kappa+g-2$ and we can argue as in case $p=3$.  Let $p_0\geq4$. Further, if $\delta-\kappa+g-|Z_3|\leq1$, then $|Z_1|\geq|Z_3|\geq\delta-\kappa+g-1$ and we can argue as in case $p=2$. So, we can assume that $\delta-\kappa+g-|Z_3|\geq2$. Since $p_0\geq4$, there are vertex disjoint paths $R_1,R_2,R_3,R_4$ in $\langle A^{\downarrow}\cup V(Q^{\downarrow}_{*}) \rangle$ connecting $P_0$ and $Q^{\downarrow}_{*}$ (otherwise, there exist a cut-set of $\langle A^{\downarrow}\cup V(Q^{\downarrow}_{*}) \rangle$ of order at most three, contradicting the definition of $A^{\downarrow}$). Let $F(R_i)\in V(P)$ and $L(R_i)\in V(Q^{\downarrow}_{*})$ $(i=1,2,3,4)$. By the choice of $w$, $|N(w_i)\cap V(P_0)|\geq\delta-\kappa+g-|Z_3|$ for each $i\in\{1,...,s\}$. In particular, for $i=s$, we have $s\geq\delta-\kappa+g-|Z_3|$. By Lemma E, in $\langle V(P_0)\rangle$ any two vertices are joined by a path of length at least $\delta-\kappa+g-|Z_3|$. Let $I_1,...,I_t$ be the segments of $Q^{\downarrow}_{*}$ having only their ends in common with $Z_1\cup Z_3\cup Z_4\cup Z_5 $, where $Z_4=\{L(R_1)\}$ and $Z_5=\{L(R_2)\}$. Form a cycle $C^{\downarrow}_{1}$ by adding to $Q^{\downarrow}_{*}$ an extra path of length $r$, where $r=\delta-\kappa+g-|Z_3|+2$, connecting $v_1$ and $w_1$. \\

\textbf{Case 2.3.1.} At least two of segments $I_1,...,I_t$ intersect $S_0$.

In this case we can apply Claim 4, which gives
$$
|V^{\downarrow}|\geq|V(Q^{\downarrow}_{*})|\geq2(|Z_1|+|Z_3|)+3r-15
$$
$$
=(3\delta-3\kappa+2)+(|Z_1|-|Z_3|)+3g+|Z_1|-11.
$$

Since $|A^{\downarrow}|\geq|V^{\downarrow}|-g+p$, we have
$$
|A^{\downarrow}|\geq(3\delta-3\kappa+2)+(|Z_1|-|Z_3|)+2g+p+|Z_1|-11>3\delta-3\kappa+2,
$$
contradicting the hypothesis.\\

\textbf{Case 2.3.2}. Only one of segments $I_1,...,I_t$, say $I_1$, intersects $S_0$.

In this case, each vertex of $S_0$ is an inner vertex for $I_1$. Form a path $L$ by deleting $I_1$ from $Q^{\downarrow}_{*}$ and adding a longest path connecting $F(I_1)$ and $L(I_1)$ and passing through $V(P)$. Then we can apply Claim 4 with respect to $L$ and can reach a contradiction as in Case 2.3.1.\\

\textbf{Case 3.} $f=2$ and $S\subseteq V^{\uparrow}$.

If $A^{\downarrow}\subseteq V^{\downarrow}$, then we are done. Let $A^{\downarrow}\not\subseteq V^{\downarrow}$ and let $P=y_1...y_p$; $P_0$, $p_0$, $w$ and $Z_1$, $Z_2$, $Z_3$ are as defined in Case 1. If $p=1$, then $|Z_1|\geq\delta-\kappa+2$ and by standard arguments, $|V^{\downarrow}|\geq2|Z_1|-1\geq2\delta-2\kappa+3$. Let $p\geq2$, implying in particular that $p_0\geq2$. If $p_0=2$, then $|Z_1|\geq|Z_2|\geq\delta-\kappa+1$ and we can argue as in Case 1.1. Let $p_0\geq3$. As above, there are vertex disjoint paths $R_1,R_2,R_3$ connecting $V^{\downarrow}$ and $V(P_0)$. Let $F(R_i)\in V(P_0)$ and $L(R_i)\in V^{\downarrow}$ $(i=1,2,3)$. If $\delta-\kappa-|Z_3|+2\leq1$, then $|Z_1|\geq|Z_3|\geq\delta-\kappa+1$ and we can argue as in the case $p_0\leq2$. Let $\delta-\kappa-|Z_3|+2\geq2$. By the choice of $w$, $|N(w^-_i)\cap V(P_0)|\geq\delta-\kappa+2-|Z_3|$ $(i=1,...,s)$. In particular, for $i=s$, we have $s\geq\delta-\kappa-|Z_3|+2$. By Lemma E, in $\langle V(P_0)\rangle$ any two vertices are joined by a path of length at least $\delta-\kappa-|Z_3|+2$. By Claim 3,
$$
|V^{\downarrow}|\geq2(|Z_1|+|Z_3|)+2r-11\geq2\delta-2\kappa+3.\qquad   \Delta 
$$

\noindent\textbf{Proof of Lemma 9}. Let $Q^{\downarrow}_{1},...,Q^{\downarrow}_{m}$ and $V^{\downarrow}_{1},...,V^{\downarrow}_{m}$, $V^{\downarrow}$ be as defined in Definition B. The existence of $Q^{\downarrow}_{1},...,Q^{\downarrow}_{m}$ follows from Lemma D. Clearly $f\geq2m$. Further, let $P=y_1...y_p$, $Z_1$, $Z_2$ and $Z_{1,i}$, $Z_{2,i}$ $(i=1,...,m)$ be as defined in proof of Lemma 8 (Case 1). \\  

\textbf{Case 1}. $f\geq3$.

Assume first that $\delta-\kappa\leq1$. By combining this with $\delta\geq2\kappa-2$, we get $f=\kappa=3$, $\delta=4$ and $m=1$. If $A^{\downarrow}\subseteq V^{\downarrow}$, then we are done. Let $A^{\downarrow}\not\subseteq V^{\downarrow}$ and let $z\in A^{\downarrow}-V^{\downarrow}$. If $z$ and $Q^{\downarrow}_1$ are connected in $\langle A^{\downarrow}\cup S\rangle$ by at most three paths having no vertex other than $z$ in common, then $G$ has a cut-set $S^{\prime}$ of order three with $S^{\prime}\subset A^{\downarrow}\cup S$ and $S^{\prime}\neq S$, contradicting the definition of $A^{\downarrow}$. Otherwise, it is easy to see that $|V^{\downarrow}|\geq7>3\delta-3\kappa+f$. So, we can assume that $\delta-\kappa\geq2$. If $p\leq3$, then we can argue as in proof of Lemma 8 (Cases 1.1-1.2). Let $p\geq4$ and let $P_0,p_0,w,Z_3$ and $Z_{3,i}$ $(i=1,...,m)$ are as defined in proof of Lemma 8 (Case 1.3). Clearly 
$$
|Z_1|\geq|Z_3|\geq|Z_2|\geq\delta-\kappa+f-p_0+1.   \eqno {(5)}
$$ 

By Claim 1 (see the proof of Lemma 8),
$$
|V^{\downarrow}|=\sum^{m}_{i=1}|V^{\downarrow}_{i}|\geq2(|Z_1|+|Z_3|)-2m\geq4(\delta-\kappa+f-p_0+1)-2m,
$$

implying that
$$
|V^{\downarrow}|\geq(3\delta-3\kappa+f)+(\delta-\kappa-2)+3(f-2m)+4(m-p_0+1)+2.   \eqno{(6)}
$$

If $p_0\leq m+1$, then $|V^{\downarrow}|\geq3\delta-3\kappa+f$ and we are done. Now let $p_0\geq m+2$. Further, if $\delta-\kappa+f-|Z_3|\leq2$, then $|Z_1|\geq|Z_3|\geq\delta-\kappa+f-2$ and by Claim 1,\\

$|V^{\downarrow}|=\sum^{m}_{i=1}|V^{\downarrow}_{i}|\geq2(|Z_1|+|Z_3|)-2m\geq4(\delta-\kappa+f-2)-2m$\\

$=(3\delta-3\kappa+f)+(\delta-\kappa-2)+3f-2m-6\geq3\delta-3\kappa+f.$\\

Now let $\delta-\kappa+f-|Z_3|\geq3$. By the choice of $w$,
$$
|N(w_i)\cap V(P_0)|\geq\delta-\kappa+f-|Z_3|\quad (i=1,...,s).
$$. 

In particular, for $i=s$, we have $s\geq\delta-\kappa+f-|Z_3|$. By Lemma E, in $\langle V(P_0)\rangle$ any two vertices are joined by a path of length at least $\delta-\kappa+f-|Z_3|$. For each $i\in \{1,...,m\}$, form a cycle,  denoted by $C^{\downarrow}_{i}$, by adding to $Q^{\downarrow}_{i}$ an extra path of length $r$, where $r=\delta-\kappa+f-|Z_3|+2\geq5$, connecting $F(Q^{\downarrow}_{i})$ and $L(Q^{\downarrow}_{i})$. Combining $\delta-\kappa+f-|Z_3|\geq3$ with (5), we get, $p_0\geq4$. \\

\textbf{Case 1.1}. $p_0=m+2$.

Substituting $p_0=m+2$ in (6), we get $|V^{\downarrow}|\geq(3\delta-3\kappa+f)+\delta-\kappa-4$. If $\delta-\kappa\geq4$, then we are done. Let $\delta-\kappa\leq3$. By the hypothesis, $\delta-\kappa\geq\kappa-2$, implying that $\kappa\leq5$, as well as $f\leq5$ and $m\leq2$. Since $m+2=p_0\geq4$, we have $m=2$ and $p_0=4$. If $f=5$, then the desired result follows from (6). If $f=4$, then it is easy to see that $|V^{\downarrow}|\geq13\geq3\delta-3\kappa+f$.\\

\textbf{Case 1.2}. $p_0\geq m+3$.

It is easy to see that $|V^{\downarrow}|\geq f+2\geq m+3$. Recalling also that $p_0\geq m+3$, we conclude that there are vertex disjoint paths $R_1,...,R_{m+3}$ in $G-(S-V^{\downarrow})$ connecting $V(P_0)$ and $V^{\downarrow}$. Since $p_0\geq m+3$, we can assume w.l.o.g. that either $L(R_i)\in V^{\downarrow}_{1}$ $(i=1,2,3,4)$ or $L(R_i)\in V^{\downarrow}_{1}$ $(i=1,2,3)$ and  $L(R_i)\in V^{\downarrow}_{2}$ $(i=4,5)$ or $L(R_i)\in V^{\downarrow}_{1}$ $(i=1,2)$, $L(R_i)\in V^{\downarrow}_{2}$ $(i=3,4)$ and $L(R_i)\in V^{\downarrow}_{3}$ $(i=5,6)$. 

\textbf{Case 1.2.1.} $L(R_i)\in V^{\downarrow}_{1}$ $(i=1,2,3,4)$.

By Claim 4 (see the proof of Lemma 8), $|V^{\downarrow}_1|\geq2(|Z_{1,1}|+|Z_{3,1}|)+3r-15$. By Claim 1, $|V^{\downarrow}_i|\geq2(|Z_{1,i}|+|Z_{3,i}|)-2$ for each $i\in \{2,...,m\}$. By summing, we get\\

$|V^{\downarrow}|=|V^{\downarrow}_1|+\sum^{m}_{i=2}|V^{\downarrow}_{i}|\geq \sum^{m}_{i=1}2(|Z_{1,i}|+|Z_{3,i}|)+3r-15-2(m-1)$\\

$= 2(|Z_1|+|Z_3|)+3(\delta-\kappa+f-|Z_3|+2)-2m-13$\\

$=(3\delta-3\kappa+f)+(|Z_1|-|Z_3|)+|Z_1|+2f-2m-7>3\delta-3\kappa+f$.\\

\textbf{Case 1.2.2.} $L(R_i)\in V^{\downarrow}_{1}$ $(i=1,2,3)$ and $L(R_i)\in V^{\downarrow}_{2}$ $(i=4,5)$.

Clearly $f\geq4$. By Claim 3, $|V^{\downarrow}_1|\geq2(|Z_{1,1}|+|Z_{3,1}|)+2r-11$. By Claim 2, $|V^{\downarrow}_2|\geq2(|Z_{1,2}|+|Z_{3,2}|)+r-7$. By Claim 1, $|V^{\downarrow}_i|\geq2(|Z_{1,i}|+|Z_{3,i}|)-2$ for each $i\in \{3,...,m\}$. By summing, we get\\

$|V^{\downarrow}|=|V^{\downarrow}_1|+|V^{\downarrow}_2|+\sum^{m}_{i=3}|V^{\downarrow}_{i}|$\\

$\geq \sum^{m}_{i=1}2(|Z_{1,i}|+|Z_{3,i}|)+3r-18-2(m-2)$\\

$= 2(|Z_1|+|Z_3|)+3(\delta-\kappa+f-|Z_3|+2)-2m-14$\\

$\geq(3\delta-3\kappa+f)+|Z_1|+2f-2m-8\geq3\delta-3\kappa+f$.\\

\textbf{Case 1.2.3.} $L(R_i)\in V^{\downarrow}_{1}$ $(i=1,2)$, $L(R_i)\in V^{\downarrow}_{2}$ $(i=3,4)$, $L(R_i)\in V^{\downarrow}_{3}$ $(i=5,6)$.

Clearly $f\geq6$. By Claim 2, $|V^{\downarrow}_i|\geq2(|Z_{1,i}|+|Z_{3,i}|)+r-7$ $(i=1,2,3)$. By Claim 1, $|V^{\downarrow}_i|\geq2(|Z_{1,i}|+|Z_{3,i}|)-2$ for each $i\in \{4,...,m\}$. By summing, we get\\

$|V^{\downarrow}|=|V^{\downarrow}_1|+|V^{\downarrow}_2|+|V^{\downarrow}_3|+\sum^{m}_{i=4}|V^{\downarrow}_{i}|$\\

$\geq \sum^{m}_{i=1}2(|Z_{1,i}|+|Z_{3,i}|)+3r-21-2(m-3)$\\

$= 2(|Z_1|+|Z_3|)+3(\delta-\kappa+f-|Z_3|+2)-2m-15$\\

$\geq(3\delta-3\kappa+f)+|Z_1|+2f-2m-9>3\delta-3\kappa+f$.\\

\textbf{Case 2}. $f=2$ and $S\not\subseteq V^{\uparrow}$.

Choose a vertex $z\in S-V^{\uparrow}$. Put $F(Q^{\downarrow}_1)=v_1$ and $L(Q^{\downarrow}_1)=w_1$. If there is a cut-vertex $x$ in $\langle A^{\downarrow}\cup\{v_1,w_1,z\}\rangle$ that separates $z$ and $Q^{\downarrow}_{1}$, then $S^{\prime }=\{S-z\}\cup \{x\}$  is a cut-set of $G$ other than $S$ with $S^{\prime}\subset A^{\downarrow}\cup S$, contradicting the definition of $A^{\downarrow}$. Otherwise, the existence of $Q^{\downarrow}_{*}$ (see Definition B) follows easily. Put $g=|V(Q^{\downarrow}_{*})\cap S|$. By the definition, $g\geq3$. Assume first that $\delta-\kappa\leq1$. Combining this with $\delta\geq2\kappa-2$, we obtain $\delta=4$ and $\kappa=3$. If $A^{\downarrow}\subseteq V(Q^{\downarrow}_{*})$, then we are done. Let $A^{\downarrow}\not\subseteq V(Q^{\downarrow}_{*})$ and choose a vertex $z\in A^{\downarrow}- V(Q^{\downarrow}_{*})$. As in Case 1, there are at least 4 paths connecting $z$ and $Q^{\downarrow}_{*}$, implying that $|V^{\downarrow}|\geq|V(Q^{\downarrow}_{*})|\geq7>3\delta-3\kappa+3$. Now let $\delta-\kappa\geq2$. Let $P=y_1...y_p;Z_1;Z_2$ be as defined in proof of Lemma 8 (Case 1). If $p\leq1$, then $A^{\downarrow}-V(Q^{\downarrow}_{*})$ is independent. Let $p\geq2$. Put $S_0=(V(Q^{\downarrow}_{*})\cap S)-\{v_1,w_1\}$.\\

\textbf{Case 2.1.} $p=2$.

Clearly $|Z_i|\geq\delta-\kappa+g-1\geq4$ $(i=1,2)$. Form a cycle $C^{\downarrow}_{1}$ by adding to $Q^{\downarrow}_{*}$ an extra path of length 3 with endvertices $v_1,w_1$. Let $I_1,...,I_t$ be the segments of $C^{\downarrow}_{1}$ having only their ends in common with $Z_1\cup Z_2$ and having at least one inner vertex in common with $S$. \\

\textbf{Case 2.1.1.} $S_0\not\subseteq V(I_i)-\{F(I_i),L(I_i)\}$ $(i=1,...,t)$.

It is easy to see that $(Z_1,Z_2)$ is a nontrivial $(C^{\downarrow}_{1},3)$-scheme. By Lemma 1,\\

$|V^{\downarrow}|\geq|V(Q^{\downarrow}_{*})|=|C^{\downarrow}_{1}|-2\geq|Z_1|+|Z_2|+|Z_1\cup Z_2|-2$\\

$\geq3(\delta-\kappa+g-1)-2=3\delta-3\kappa+3g-5>3\delta-3\kappa+3$.\\

\textbf{Case 2.1.2.} $S_0\subseteq V(I_i)-\{F(I_i),L(I_i)\}$ for some $i\in \{1,...,t\}$, say $i=1$.

Let $I_1=x\overrightarrow{I_1}y$. Denote by $R$ a longest path connecting $x$ and $y$ and passing through $V(P)$. Clearly $|R|=3$ if $x$ and $y$ belong to different $Z_1,Z_2$ and $|R|=2$, otherwise. Form a new cycle $C^{\downarrow}_{2}$ by deleting $I_1$ from $C^{\downarrow}_{1}$ and adding $R$. Clearly $(Z_1,Z_2)$ is a nontrivial $(C^{\downarrow}_{2},3)$-scheme and as in Case 2.1.1, $|V^{\downarrow}|>3\delta-3\kappa+3$.\\

\textbf{Case 2.2.} $p=3$.

Clearly $|Z_i|\geq\delta-\kappa+g-2\geq3$ $(i=1,2)$. Form a cycle $C^{\downarrow}_{1}$ by adding to $Q^{\downarrow}_{*}$ an extra path of length 4 with endvertices $v_1,w_1$. Let $I_1,...,I_t$ and $S_0$ be as defined in Case 2.1.\\

\textbf{Case 2.2.1.} $S_0\not\subseteq V(I_i)-\{F(I_i),L(I_i)\}$ $(i=1,...,t)$.

It is easy to see that $(Z_1,Z_2)$ is a nontrivial $(C^{\downarrow}_{1},4)$-scheme. By Lemma 1,\\

$|V^{\downarrow}|\geq|V(Q^{\downarrow}_{*})|=|C^{\downarrow}_{1}|-3\geq2(|Z_1|+|Z_2|)-3$\\

$\geq4(\delta-\kappa+g-2)-3=3\delta-3\kappa+3+(\delta-\kappa)+4g-14\geq3\delta-3\kappa+3$.\\

\textbf{Case 2.2.2.} $S_0\subseteq V(I_i)-\{F(I_i),L(I_i)\}$ for some $i\in \{1,...,t\}$, say $i=1$.

Let $I_1=x\overrightarrow{I_1}y$. Denote by $R$ a longest path connecting $x$ and $y$ and passing through $V(P)$. Clearly $|R|=4$ if $x$ and $y$ belong to different $Z_1,Z_2$ and $|R|\geq2$, otherwise. Form a new cycle $C^{\downarrow}_{2}$ by deleting $I_1$ from $C^{\downarrow}_{1}$ and adding $R$. Clearly $(Z_1,Z_2)$ is a nontrivial $(C^{\downarrow}_{2},4)$-scheme and as in Case 2.2.1, $|V^{\downarrow}|\geq3\delta-3\kappa+3$.\\

\textbf{Case 2.3.} $p\geq4$.

Let $P_0,p_0,w$ and $Z_3$ are as defined in proof of Lemma 8 (Case 1.3). By the definition, $|Z_1|\geq|Z_3|\geq|Z_2|$ and $p_0\geq2$. If $p_0=2$, then $|Z_1|\geq|Z_2|\geq\delta-\kappa+g-1$ and we can argue as in case $p=2$. If $p_0=3$, then $|Z_1|\geq|Z_2|\geq\delta-\kappa+g-2$ and we can argue as in case $p=3$.  Let $p_0\geq4$. Further, if $\delta-\kappa+g-|Z_3|\leq1$, then $|Z_1|\geq|Z_3|\geq\delta-\kappa+g-1$ and we can argue as in case $p=2$. Let $\delta-\kappa+g-|Z_3|\geq2$. Since $p_0\geq4$, there are vertex disjoint paths $R_1,R_2,R_3,R_4$ in $\langle A^{\downarrow}\cup V(Q^{\downarrow}_{*}) \rangle$ connecting $P_0$ and $Q^{\downarrow}_{*}$ (otherwise, there exist a cut-set of $\langle A^{\downarrow}\cup V(Q^{\downarrow}_{*}) \rangle$ of order 3 contradicting the definition of $A^{\downarrow}$). Let $F(R_i)\in V(P)$ and $L(R_i)\in V(Q^{\downarrow}_{*})$ $(i=1,2,3,4)$. Clearly, $|N(w_i)\cap V(P_0)|\geq\delta-\kappa+g-|Z_3|$ for each $i\in\{1,...,s\}$. In particular, for $i=s$, we have $s\geq\delta-\kappa+g-|Z_3|$. By Lemma E, in $\langle V(P_0)\rangle$ any two vertices are joined by a path of length at least $\delta-\kappa+g-|Z_3|$. Assume first that $|Z_3|\leq3$ and assume w.l.o.g. that $L(R_1),L(R_2),L(R_3),L(R_4)$ occur on $v_1\overrightarrow{Q}^{\downarrow}_{*}w_1$ in a consecutive order. Consider the segments $J_i=L(R_i)\overrightarrow{Q}^{\downarrow}_{*}L(R_{i+1})$, $i=1,2,3$. If at least two of $J_1,J_2,J_3$ intersect $S_0$, then by the definition of $Q^{\downarrow}_{*}$, $|J_i|\geq\delta-\kappa+g-|Z_3|+2\geq4$ $(i=1,2,3)$ and hence,\\

$|V^{\downarrow}|\geq|V(Q^{\downarrow}_{*})|\geq|J_1|+|J_2|+|J_3|+1\geq3(\delta-\kappa+g-|Z_3|+2)+1$\\

$=(3\delta-3\kappa+3)+3g-3|Z_3|+4>3\delta-3\kappa+3$.\\

Now let only one of $J_1,J_2,J_3$, say $J_1$, intersects $S_0$. Form a new path $L$ by deleting $J_1$ from $Q^{\downarrow}_{*}$ and adding a path of length $\delta-\kappa+g-|Z_3|+2$, connecting $L(R_1)$ and $L(R_2)$ and passing through $V(P_0)$. Then as above, $|V^{\downarrow}|>|L|\geq3\delta-3\kappa+3$. Now let $|Z_3|\geq4$, implying that $|Z_1|\geq|Z_3|\geq4$. Assume w.l.o.g. that $R_1\cup R_2$ does not intersect $y_1\overrightarrow{P}w^{-}_{1}$ (see the proof of Lemma 8, Case 1.3) and does not contain $w$. Let $I_1,...,I_t$ be the segments of $Q^{\downarrow}_{*}$ having only their ends in common with $Z_1\cup Z_3\cup Z_4\cup Z_5 $, where $Z_4=\{L(R_1)\}$ and $Z_5=\{L(R_2)\}$. Form a cycle $C^{\downarrow}_{1}$ by adding to $Q^{\downarrow}_{*}$ an extra path of length $r$, where $r=\delta-\kappa+g-|Z_3|+2$, connecting $v_1$ and $w_1$. \\

\textbf{Case 2.3.1.} At least two of segments $I_1,...,I_t$ intersect $S_0$.

It is easy to see that $(Z_1,Z_3,Z_4,Z_5)$ is a nontrivial $(C^{\downarrow}_{1},r$-scheme. By $(c1)$ ( see Lemma 3), \\

$|V^{\downarrow}|\geq|V(Q^{\downarrow}_{*})|=|C^{\downarrow}_{1}|-r+1\geq2(|Z_1|+|Z_3|)+3r-15$\\

$=(3\delta-3\kappa+3)+(|Z_1|-|Z_3|)+3g+|Z_1|-12>3\delta-3\kappa+3$.\\

\textbf{Case 2.3.2}. Only one of segments $I_1,...,I_t$, say $I_1$, intersects $S_0$.

In this case, each vertex of $S_0$ is an inner vertex for $I_1$. Let $I_1=xI_1y$. Form a cycle $C^{\downarrow}_{2}$ by deleting $I_1$ from $C^{\downarrow}_{1}$ and adding a longest path connecting $x$ and $y$ and passing through $V(P)$. It is easy to see that $(Z_1,Z_3,Z_4,Z_5)$ is a nontrivial $(C^{\downarrow}_{1},r)$-scheme. Then we can argue as in Case 2.2.1.\\

\textbf{Case 3.} $f=2$ and $S\subseteq V^{\uparrow}$.

If $A^{\downarrow}\subseteq V^{\downarrow}$, then $|V^{\downarrow}|\geq|A^{\downarrow}|+2\geq3\delta-3\kappa+4>2\delta-2\kappa+3$. Let $A^{\downarrow}\not\subseteq V^{\downarrow}$. Then we can argue exactly as in proof of Lemma 8 (Case 3). \qquad $\Delta$\\

\noindent\textbf{Proof of Lemma 10}. To prove $(f1)$, choose a vertex $z\in A^{\uparrow}-V^{\uparrow}$. Clearly $N(z)\subseteq V^{\uparrow}\cup S$. Put\\
 
$Z_i=N(z)\cap V^{\uparrow}_i \quad (i=1,...,m), \quad Z=\cup^{m}_{i=1}Z_{i},$\\

$W=N(z)\cap \{v_1,...,v_m,w_1,...,w_m\}.$\\

Since $\{Q^{\uparrow}_{1},...,Q^{\uparrow}_{m}\}$ is $(*1)$-extreme, we have $|W|\leq2$ and $|Z|\geq\delta-1.$ If $|W|=0$, then by standard arguments, $|V^{\uparrow}_{i}|\geq2|Z_i|+1$ $(i=1,...,m)$ and
$$
|V^{\uparrow}|=\sum^{m}_{i=1}|V^{\uparrow}_{i}|\geq\sum^{m}_{i=1}(2|Z_i|+1)=2|Z|+m\geq2\delta+m-2.
$$

Further, let $|W|=1$. Assume w.l.o.g. that $W=\{v_1\}$. By $(*1)$, $N(z)\subseteq V^{\uparrow}$, i.e. $|Z|\geq\delta$. Then clearly $|V^{\uparrow}_{1}|\geq2|Z_1|$ and $|V^{\uparrow}_{i}|\geq2|Z_i|+1$ $(i=2,...,m)$, implying that\\

$|V^{\uparrow}|= |V^{\uparrow}_{1}|+\sum^{m}_{i=2}|V^{\uparrow}_{i}|\geq 2|Z_1|+\sum^{m}_{i=2}(2|Z_i|+1)$\\

$\geq2|Z|+m-1>2\delta+m-2.$\\

Finally, we assume that $|W|=2$. By $(*1)$, $W\subseteq V^{\uparrow}_{i}$ for some $i\in \{1,...,m\}$, say $i=1$, and $N(z)\subseteq V^{\uparrow}$ implying that $|Z|\geq\delta$. Then clearly $|V^{\uparrow}_{1}|\geq2|Z_1|-1$ and $|V^{\uparrow}_{i}|\geq2|Z_i|+1$ $(i=2,...,m)$, whence\\

$|V^{\uparrow}|= |V^{\uparrow}_{1}|+\sum^{m}_{i=2}|V^{\uparrow}_{i}|\geq 2|Z_1|-1+\sum^{m}_{i=2}(2|Z_i|+1)$\\

$=2|Z|+m-2\geq2\delta+m-2.$\\

To prove $(f2)$, choose a vertex $z\in A^{\downarrow}-V^{\downarrow}$ and put 
$$
Z_i=N(z)\cap V^{\downarrow}_i \quad (i=1,...,m), \quad Z=\cup^{m}_{i=1}Z_{i},
$$

Clearly $|Z|\geq\delta-(\kappa-f)$. Since $\{Q^{\downarrow}_{1},...,Q^{\downarrow}_{m}\}$ is $(*2)$-extreme, by standard arguments, $|V^{\downarrow}_{i}|\geq2|Z_i|-1$ $(i=1,...,m)$, implying that 
$$
|V^{\downarrow}|=\sum^{m}_{i=1}|V^{\downarrow}_{i}|\geq\sum^{m}_{i=1}(2|Z_i|-1)=2|Z|-m\geq2\delta-2\kappa+2f-m.  
$$

To prove $(f3)$, choose a vertex $z\in A^{\downarrow}-V(Q^{\downarrow}_{*})$. Put $Z=N(z)\cap V(Q^{\downarrow}_{*})$ and $g=|V(Q^{\downarrow}_{*})\cap S|$. Clearly $|Z|\geq\delta-\kappa+g\geq\delta-\kappa+3$. By standard arguments, $|V(Q^{\downarrow}_{*})|\geq2|Z|-1\geq2\delta-2\kappa+5$. \qquad $\Delta$

\section{Proofs of theorems}

\textbf{Proof of Theorem 2}. If $G$ has a dominating cycle, then we are done. Otherwise, by Lemma M, $c\geq 3\delta-3$. If $\delta\leq 2\kappa-3$, then $c\geq 3\delta-3\geq4\delta-2\kappa$ and again we are done. So, we can assume that $\delta\geq 2\kappa-2>3\kappa/2-1$. Let $Q^{\uparrow}_{1},...,Q^{\uparrow}_{m}$; $Q^{\downarrow}_{1},...,Q^{\downarrow}_{m}$ and $C$ be as defined in Definitions A-C. The existence of $Q^{\downarrow}_{1},...,Q^{\downarrow}_{m}$ follows from Lemma D. Put $f=|V^{\downarrow}\cap S|$. Assume w.l.o.g. that $F(Q^{\uparrow}_{i})=v_i$ and $L(Q^{\uparrow}_{i})=w_i$ $(i=1,...,m)$. Form a cycle $C^{\uparrow}$ consisting of $Q^{\uparrow}_{1},...,Q^{\uparrow}_{m}$ and extra edges $w_1v_2, w_2v_3,...,w_{m-1}v_m,w_mv_1$. By Lemma 5, $A^{\uparrow}-V^{\uparrow}$ is independent.\\ 

\textbf{Case 1}. $|A^\downarrow|\leq2\delta-2\kappa+1.$

By Lemma 7, $A^{\downarrow}\subseteq V^{\downarrow}$. By $(*1)$, $S-V(C)$ is independent. If $V(G-C)$ is independent, then clearly $C$ is a dominating cycle. Otherwise, each edge of $G-C$ has one end in $A^{\uparrow}-V^{\uparrow}$ and the other in $S-V(C)$. Let $xy\in E(G)$, where $x\in A^{\uparrow}-V^{\uparrow}$ and $y\in S-V(C)$. By $(*1)$, $N(x)\subseteq V^{\uparrow}\cup\{y\}$.\\

\textbf{Case 1.1.} $N(y)\not\subseteq V(C)\cup \{x\}$.

Let $yz\in E(G)$ for a vertex $z\in A^{\uparrow}-V^{\uparrow}$ other than $x$. By $(*1)$ and the fact that $A^{\uparrow}-V^{\uparrow}$ is independent, we have $N(z)\subseteq V^{\uparrow}\cup\{y\}$. Set $Z_1=N(x)\cap V^{\uparrow}$ and $Z_2=N(z)\cap V^{\uparrow}$. Clearly $|Z_i|\geq\delta-1$ $(i=1,2)$. Basing on $(*1)$, it is easy to see that $(Z_1,Z_2)$ is a nontrivial $(C^{\uparrow},4)$-scheme. By $(a1)$ (see Lemma 1), $|C|>|C^{\uparrow}|\geq4(\delta-1)>4\delta-2\kappa$.\\

\textbf{Case 1.2.} $N(y)\subseteq V(C)\cup \{x\}$.

\textbf{Case 1.2.1.} $|N(y)\cap V^{\downarrow}|\geq f-m+1.$  

Set $Z_{i}=N(y)\cap V^{\downarrow}_{i}$ $(i=1,...,m)$ and $Z=\cup^{m}_{i=1}Z_{i}$. By the hypothesis, $|Z|=|N(y)\cap V^{\downarrow}|\geq f-m+1$. By $(*1)$, $y$ is not adjacent to $F(Q^{\downarrow}_{i})$ and $L(Q^{\downarrow}_{i})$ $(i=1,...,m)$. Since $\{Q^{\downarrow}_{1},...,Q^{\downarrow}_{m}\}$ is $(*4)$-extreme, by standard arguments, $|V^{\downarrow}_{i}\cap S|\geq|Z_{i}|+1$ $(i=1,...,m)$ and
$$ 
f=|V^{\downarrow}\cap S|=\sum^{m}_{i=1}|V^{\downarrow}_{i}\cap S|\geq \sum^{m}_{i=1}(|Z_i|+1)\geq|Z|+m\geq f+1,
$$
a contradiction. \\

\textbf{Case 1.2.2.} $|N(y)\cap V^{\downarrow}|\leq f-m\leq f-1$.

Since $N(y)\subseteq V(C)\cup \{x\}$, we have $|N(y)\cap V^{\uparrow}|\geq \delta-|N(y)\cap V^{\downarrow}|-1\geq\delta-f.$ Set $Z_1=N(x)\cap V^{\uparrow}$ and $Z_2=N(y)\cap V^{\uparrow}$. Clearly $|Z_1|\geq\delta-1$ and $|Z_2|\geq\delta-f$. By $(*1)$, it is easy to see that $(Z_1,Z_2)$ is a nontrivial $(C^{\uparrow},3)$-scheme. By Lemma 1, $|V^{\uparrow}|=|C^{\uparrow}|\geq2(\delta-1)+\delta-f=3\delta-f-2$. Observing that $|A^{\downarrow}|\geq\delta-\kappa+1$, we obtain
$$
|C|\geq|V^{\uparrow}|+|A^{\downarrow}|\geq (4\delta-2\kappa)+\kappa-f-1\geq4\delta-2\kappa.
$$.

\textbf{Case 2}. $|A^\downarrow|\geq2\delta-2\kappa+2.$

\textbf{Case 2.1.} $f\geq3$.

By Lemma 8 and $(*1)$, $A^{\downarrow}-V^{\downarrow}$ and $S-V(C)$ both are independent sets. If $A^{\downarrow}\subseteq V^{\downarrow}$, then we can argue exactly as in Case 1. Let $A^{\downarrow}\not\subseteq V^{\downarrow}$.
By $(f2)$ (see Lemma 10), $|V^{\downarrow}|\geq2\delta-2\kappa+2f-m.$ If $A^{\uparrow}\not\subseteq V^{\uparrow}$, then by $(f1)$ (see Lemma 10), $|V^{\uparrow}|\geq2\delta+m-2$, whence
$$
|C|\geq|V^{\uparrow}|+|V^{\downarrow}|-2m\geq(4\delta -2\kappa)+2f-2m-2>4\delta-2\kappa.
$$

Now let $A^{\uparrow}\subseteq V^{\uparrow}$. Since $S-V(C)$ is independent, either $C$ is a dominating cycle (then we are done) or each edge in $G-C$ has one end in $S-V(C)$ and the other in $A^{\downarrow}-V^{\downarrow}$.\\

\textbf{Case 2.1.1.} There is a path $xyz$ in $G-C$ with $x,z\in A^{\downarrow}$ and $y\in S$.  

Put $Z_1=N(x)\cap V^{\downarrow}$, $Z_2=N(z)\cap V^{\downarrow}$ and
$$
Z_{1,i}=N(x)\cap V^{\downarrow}_{i}, \quad Z_{2,i}=N(z)\cap V^{\downarrow}_{i}\quad (i=1,...,m).
$$

Clearly $|Z_i|\geq\delta-\kappa+f$ $(i=1,2)$. For each $i\in \{1,...,m\}$, form a cycle $C^{\downarrow}_{i}$ by adding to $Q^{\downarrow}_{i}$ an extra path of length 4 with endvertices $F(Q^{\downarrow}_{i})$ and $L(Q^{\downarrow}_{i})$. If $(Z_{1,i},Z_{2,i})$ is a nontrivial $(Q^{\downarrow}_{i},4)$-scheme, then by $(a1)$ (see Lemma 1),  $|V^{\downarrow}_{i}|=|C^{\downarrow}_{i}|-3\geq2(|Z_{1,i}|+|Z_{2,i}|)-3$. Otherwise, it can be checked easily. By summing, we get\\

$|V^{\downarrow}|=\sum^{m}_{i=1}|V^{\downarrow}_{i}|\geq\sum^{m}_{i=1}(2(|Z_{1,i}|+|Z_{2,i}|)-3)=2(|Z_1|+|Z_2|)-3m$\\

$\geq 4(\delta-\kappa+f)-3m \geq 3\delta-3\kappa+f+2,$\\

\noindent whence $|A^{\downarrow}|\geq|V^{\downarrow}|-f > 3\delta-3\kappa+2$, contrary to the hypothesis.\\

\textbf{Case 2.1.2.} There is a path $xyz$ in $G-C$ with $x,z\in S$ and $y\in A^{\downarrow}$.

Clearly $\kappa\geq4$. If $xy^{\prime}\in E(G)$ or $zy^{\prime}\in E(G)$ for some $y^{\prime}\in A^{\downarrow}-V^{\downarrow}$ other than $y$, then we can argue as in Case 2.1.1. So, we can assume that $N(x)\subseteq V(C)\cup \{y\}$ and $N(z)\subseteq V(C)\cup \{y\}$. Put $Z_1=N(x)\cap V(C)$ and $Z_2=N(z)\cap V(C)$. \\

\textbf{Case 2.1.2.1.} $m=1$.

Assume w.l.o.g. that $Q^{\uparrow}_{1}=v_1\overrightarrow{C}w_1$. Let $I=\xi\overrightarrow{C}\eta$ be a segment on $C$ having only $\xi,\eta$ in common with $Z_1\cup Z_2$. By $(*1)$, $\{\xi,\eta\}\cap \{v_1,w_1\}=\emptyset.$ If $\xi,\eta\in Z_1$ or $\xi,\eta\in Z_2$, then  by $(*1)$ and $(*2)$, $|I|\geq2$. Let $\xi\in Z_1$ and $\eta\in Z_2$. If $\xi,\eta\in V^{\downarrow}$, then by $(*2)$, $|I|\geq4$. Let $\xi,\eta\in V^{\uparrow}$. If $|I|\leq2$, then the collection of paths $v_1\overrightarrow{C}\xi x$, $z\eta\overrightarrow{C}w_1$ contradicts $(*1)$. If $|I|=3$, then the cycle $v_1\overrightarrow{C}\xi xyz\eta\overrightarrow{C}v_1$ contradicts $(*4)$. So, $|I|\geq4$ if $\xi,\eta\in V^{\uparrow}$. As for the cases $\xi\in V^{\uparrow}$, $\eta\in V^{\downarrow}$ or $\xi\in V^{\downarrow}$, $\eta\in V^{\uparrow}$, there exists exactly one segment $I_1=\xi_1\overrightarrow{C}\eta_1$ with $\xi_1\in V^{\uparrow}$, $\eta_1\in V^{\downarrow}$, $w_1\in I_1$ and exactly one segment $I_2=\xi_2\overrightarrow{C}\eta_2$ with $\xi_2\in V^{\downarrow}$, $\eta_2\in V^{\uparrow}$, $v_1\in I_2$, each of length at least 2. Form a new cycle $C^{\prime}$ by deleting $I_1$ and $I_2$ from $C$ and adding appropriate extra paths each of length 4. Clearly $(Z_1,Z_2)$ is a nontrivial $(C^{\prime},4)$-scheme. By $(a1)$ (see Lemma 1),
$$
|C|\geq |C^{\prime}|-4\geq 2(|Z_1|+|Z_2|)-4\geq 4(\delta-1)-4\geq 4\delta-2\kappa.
$$

\textbf{Case 2.1.2.2.} $m\geq2$.

By arguing as in Case 2.1.2.1, we obtain
$$
|C|\geq |C^{\prime}|-4m\geq 4(\delta-1)-4m=(4\delta-2\kappa)+2(\kappa-2m-2).
$$

Since $\kappa\geq f+2\geq2m+2$, we have $|C|\geq4\delta-2\kappa$.\\
 
\textbf{Case 2.1.3.} $xy\in E(G)$ for some $x\in S-V(C)$ and $y\in A^{\downarrow}-V^{\downarrow}.$

Form a new cycle $C^{\downarrow}$ by deleting $Q^{\uparrow}_{1},...,Q^{\uparrow}_{m}$ from $C$ and adding extra paths $v_iz_iw_i$ $(i=1,...,m)$. Put $V_{0}=\{v_1,...,v_m,w_1,...,w_m\}$ and 
$$
Z_0=N(y)\cap V_0, \quad Z_1=N(x)\cap V^{\uparrow}, \quad Z_2=N(x)\cap V^{\downarrow},
$$
$$
Z_3=N(y)\cap (V^{\uparrow}-V_0), \quad Z_4=N(y)\cap (V^{\downarrow}-V_0).
$$

If $N(x)\not\subseteq V(C)\cup\{y\}$ or $N(y)\not\subseteq V(C)\cup\{x\}$, then we can argue as in Cases 2.1.1-2.1.2. Let $N(x)\subseteq V(C)\cup\{y\}$ and $N(y)\subseteq V(C)\cup\{x\}$. Abbreviate $|Z_0|=t$, $|Z_2|=d$ and $|Z_3|=h$. Clearly $|Z_1|\geq\delta-d-1$ and $|Z_4|\geq \delta-h-t-1.$ If $d=0$, then $|Z_1|\geq\delta-1$ and arguing exactly as in proof of $(f1)$ (see Lemma 10), we get $|V^{\uparrow}|\geq2\delta+m-2$. By $(f2)$ (see Lemma 10), $|V^{\downarrow}|\geq2\delta-2\kappa+2f-m$. Hence,
$$
|C|\geq |V^{\uparrow}|+|V^{\downarrow}|-2m\geq(4\delta-2\kappa)+2f-2m-2> 4\delta-2\kappa.
$$

Let $d\geq1$. If $h=0$, then it is not hard to see that $(Z_2,Z_0\cup Z_4)$ is a nontrivial $(C^{\downarrow},3)$-scheme. By Lemma 1, 
$$
|V^{\downarrow}|\geq |C^{\downarrow}|-m\geq 2|Z_0\cup Z_4|+|Z_2|-m\geq2(\delta-1)+d-m\geq2\delta-1-m.
$$

Further, $|A^{\downarrow}|\geq |V^{\downarrow}|-f \geq 2\delta-1-f-m$, implying that $|A^{\uparrow}|\geq2\delta-1-f-m$ and
$$
|C|\geq |A^{\uparrow}|+|V^{\downarrow}|\geq (4\delta-2\kappa)+2\kappa-f-2m-2.
$$

Observing that $\kappa\geq f+1$ and $f\geq2m$, we get $|C|\geq 4\delta-2\kappa.$ Now let $h\geq1$. Basing on $(*1)$ and $(*4)$, it is not hard to see that $(Z_1,Z_0\cup Z_3)$ is a nontrivial $(C^{\uparrow},3)$-scheme. Since $|Z_1|\geq\delta-d-1$ and $|Z_0\cup Z_3|=t+h$, we have by Lemma 1,
$$
|V^{\uparrow}|=|C^{\uparrow}|\geq 2|Z_1|+|Z_0\cup Z_3|\geq2\delta-2d-2+t+h.
$$

On the other hand, by $(*2)$, $(Z_0\cup Z_4,Z_2)$ is a nontrivial $(C^{\downarrow},3)$-scheme. By Lemma 1, 
$$
|V^{\downarrow}|=|C^{\downarrow}|-m\geq 2|Z_0\cup Z_4|+|Z_2|-m\geq 2\delta-2h+d-m-2.
$$

Note that for the special case $m=1$ and $t=0$, by the same arguments, we can obtain a slightly better estimation, namely $|V^{\downarrow}|\geq2\delta-2h+d-2$. So,
$$
|C|\geq |V^{\uparrow}|+|V^{\downarrow}|-2m\geq(4\delta-2\kappa)+2\kappa-d-h-3m+t-4.
$$

If $d\leq 2\kappa-h-3m+t-4$, then clearly $c\geq 4\delta-2\kappa$. Let $d\geq 2\kappa-h-3m+t-3$. Therefore,
$$
|V^{\downarrow}|\geq2\delta-2h+d-m-2\geq 2\delta+2\kappa-3h-4m+t-5
$$ 
whence
$$
|A^{\uparrow}|\geq |A^{\downarrow}|\geq |V^{\downarrow}|-f+1\geq2\delta+2\kappa-3h-4m+t-f-4.
$$

Further, 
$$
|C|\geq |A^{\uparrow}|+|V^{\downarrow}|\geq(4\delta-2\kappa)+6\kappa-6h-8m+2t-f-9.
$$
Observing that $\kappa\geq f+h+1$, we obtain $6\kappa\geq6(f+h+1)\geq6h+f+10m+6$, whence $|C|\geq(4\delta-2\kappa)+2m+2t-3$. If $m+t\geq2$, then clearly $|C|\geq4\delta-2\kappa$. Otherwise, we have $m=1$ and $t=0$. Recalling that in this special case, $|V^{\downarrow}|\geq 2\delta-2h+d-2$, we can obtain the desired result by the same arguments used above.\\

\textbf{Case 2.2.} $f=2$ and $S\not\subseteq V^{\uparrow}$.

By Lemma 8, $A^{\downarrow}-V(Q^{\downarrow}_{*})$ is independent. If $A^{\uparrow}\subseteq V^{\uparrow}$ and $A^{\downarrow}\subseteq V(Q^{\downarrow}_{*})$, then clearly $C$ is a dominating cycle. Next, if $A^{\uparrow}\subseteq V^{\uparrow}$ and $A^{\downarrow}\not\subseteq V(Q^{\downarrow}_{*})$, then we can argue exactly as in Case 2.1.3. Let $A^{\uparrow}\not\subseteq V^{\uparrow}$ and $A^{\downarrow}\subseteq V(Q^{\downarrow}_{*})$. By $(f1)$ (see Lemma 10), $|V^{\uparrow}|\geq2\delta-1$. Recalling that $|A^{\downarrow}|\geq 2\delta-2\kappa+2$, we get
$$
|C|\geq |V^{\uparrow}|+|A^{\downarrow}|\geq (4\delta-2\kappa)+1>4\delta-2\kappa.
$$

If $A^{\uparrow}\not\subseteq V^{\uparrow}$ and $A^{\downarrow}\not\subseteq V(Q^{\downarrow}_{*})$, then by $(f1)$ and $(f3)$ (see Lemma 10), $|V^{\uparrow}|\geq2\delta-1$ and $|V^{\downarrow}|\geq2\delta-2\kappa+5$, whence
$$
|C|\geq |V^{\uparrow}|+|V^{\downarrow}|-2\geq 4\delta-2\kappa+2>4\delta-2\kappa.
$$

\textbf{Case 2.3.} $f=2$ and $S\subseteq V^{\uparrow}.$

By Lemma 8, either $A^{\downarrow}-V^{\downarrow}=\emptyset$ or $|V^{\downarrow}|\geq2\delta-2\kappa+3$. If $A^{\downarrow}-V^{\downarrow}=\emptyset$, then $C$ is a dominating cycle, since $A^{\uparrow}-V^{\uparrow}$ is independent (by Lemma 5). So, we can assume that $|V^{\downarrow}|\geq2\delta-2\kappa+3$. If $A^{\uparrow}\not\subseteq V^{\uparrow}$, then by $(f1)$ (see Lemma 10), $|V^{\uparrow}|\geq2\delta-1$ and $|C|\geq |V^{\uparrow}|+|V^{\downarrow}|-2\geq 4\delta-2\kappa$. Now let $A^{\uparrow}\subseteq V^{\uparrow}$. If $|A^{\downarrow}|\leq3\delta-3\kappa$, then by Lemma F, $A^{\downarrow}-V^{\downarrow}$ is independent and $C$ is a dominating cycle. Let $|A^{\downarrow}|=3\delta-3\kappa+1\geq2\delta-\kappa-1$, implying that $|A^{\uparrow}|\geq2\delta-\kappa-1$. Then $|V^{\uparrow}|\geq |A^{\uparrow}|+\kappa\geq 2\delta-1$ and 
$$
|C|\geq |V^{\uparrow}|+|V^{\downarrow}|-2\geq 4\delta-2\kappa. \qquad \Delta
$$

\noindent\textbf{Proof of Theorem 3}. If $\delta\leq2\kappa-3$, then we are done by Lemma 4. Let $\delta\geq2\kappa-2$ and let $Q^{\uparrow}_{1},...,Q^{\uparrow}_{m}$; $Q^{\downarrow}_{1},...,Q^{\downarrow}_{m}$ and $C$ be as defined in Definitions A-C. The existence of $Q^{\downarrow}_{1},...,Q^{\downarrow}_{m}$ follows from Lemma D. Assume w.l.o.g. that $F(Q^{\uparrow}_{i})=v_i$ and $L(Q^{\uparrow}_{i})=w_i$ $(i=1,...,m)$. Form a cycle $C^{\uparrow}$ consisting of $Q^{\uparrow}_{1},...,Q^{\uparrow}_{m}$ and extra edges $w_1v_2, w_2v_3,...,w_{m-1}v_m,w_mv_1$. By Lemma 5, $A^{\uparrow}-V^{\uparrow}$ is independent. Put $f=|V^{\downarrow}\cap S|$. By the hypothesis, $|A^{\uparrow}|\geq|A^{\downarrow}|\geq3\delta-3\kappa+2\geq2\delta-\kappa$. \\

\textbf{Case 1}. $f\geq3$.

By Lemma 9, either $A^{\downarrow}-V^{\downarrow}$ is independent or $|V^{\downarrow}|\geq3\delta-3\kappa+f$. If $A^{\downarrow}-V^{\downarrow}$ is independent, then we can argue as in proof of Theorem 2 (Case 2.1). Let $|V^{\downarrow}|\geq3\delta-3\kappa+f\geq2\delta-2\kappa+1$. If $A^{\uparrow}\subseteq V^{\uparrow}$, then
$$
|C|\geq |A^{\uparrow}|+|V^{\downarrow}|\geq(2\delta-\kappa)+(2\delta-\kappa+1)>4\delta-2\kappa.
$$

Otherwise, by $(f1)$ (see Lemma 10), $|V^{\uparrow}|\geq2\delta+m-2$, whence 
$$
|C|\geq |V^{\uparrow}|+|V^{\downarrow}|-2m\geq(4\delta-2\kappa)+\delta-\kappa+f-m-2\geq4\delta-2\kappa.
$$

\textbf{Case 2}. $f=2$ and $S\not\subseteq V^{\uparrow}$.

By Lemma 9, either $A^{\downarrow}-V(Q^{\downarrow}_{*})$ is independent or $|V^{\downarrow}|\geq3\delta-3\kappa+3$. If $A^{\downarrow}-V(Q^{\downarrow}_{*})$ is independent, then we can argue as in proof of Theorem 2 (Case 2.2). If $|V^{\downarrow}|\geq3\delta-3\kappa+3$, then we can argue as in case $f\geq3$. \\

\textbf{Case 3}. $f=2$ and $S\subseteq V^{\uparrow}$.

By Lemma 9, $|V^{\downarrow}|\geq2\delta-2\kappa+3$. If $A^{\uparrow}\subseteq V^{\uparrow}$, then
$$
|C|\geq |A^{\uparrow}|+|V^{\downarrow}|+\kappa-2\geq(4\delta-2\kappa)+1>4\delta-2\kappa.
$$

Otherwise, by $(f1)$ (see Lemma 10),  $|V^{\uparrow}|\geq2\delta-1$, whence
$$
|C|\geq |V^{\uparrow}|+|V^{\downarrow}|-2\geq4\delta-2\kappa. \qquad   \Delta
$$

\noindent\textbf{Proof of Theorem 4}. If $\delta\leq2\kappa-3$, then we are done by Lemma 4. Let $\delta\geq2\kappa-2$ and let $Q^{\uparrow}_{1},...,Q^{\uparrow}_{m}$; $Q^{\downarrow}_{1},...,Q^{\downarrow}_{m}$ and $C$ be as defined in Definitions A-C. The existence of $Q^{\downarrow}_{1},...,Q^{\downarrow}_{m}$ follows from Lemma D. Assume w.l.o.g. that $F(Q^{\uparrow}_{i})=v_i$ and $L(Q^{\uparrow}_{i})=w_i$ $(i=1,...,m)$. Form a cycle $C^{\uparrow}$ consisting of $Q^{\uparrow}_{1},...,Q^{\uparrow}_{m}$ and extra edges $w_1v_2, w_2v_3,...,w_{m-1}v_m,w_mv_1$. By Lemmas 6, either $A^{\uparrow}-V^{\uparrow}$ is independent or $|V^{\uparrow}|\geq3\delta-5$. If $A^{\uparrow}-V^{\uparrow}$ is independent, then we can argue as in proof of Theorem 2. Let $|V^{\uparrow}|\geq3\delta-5$. \\

\textbf{Case 1.} $f\geq3.$

By Lemma 8, $A^{\downarrow}-V^{\downarrow}$ is independent. If $A^{\downarrow}\subseteq V^{\downarrow}$, then
$$
|C|\geq |V^{\uparrow}|+|A^{\downarrow}|\geq (3\delta-5)+(\delta-\kappa+1)=4\delta-\kappa-4\geq4\delta-2\kappa. 
$$

Otherwise, by $(f2)$ (see Lemma 10), $|V^{\downarrow}|\geq 2\delta-2\kappa+2f-m$, whence
$$
|C|\geq |V^{\uparrow}|+|V^{\downarrow}|-2m\geq (4\delta-2\kappa)+\delta+2f-3m-5\geq4\delta-2\kappa. 
$$

\textbf{Case 2.} $f=2$ and $S\not\subseteq V^{\uparrow}.$

By Lemma 8, $A^{\downarrow}-V(Q^{\downarrow}_*)$ is independent. If $A^{\downarrow}\subseteq V(Q^{\downarrow}_*)$, then we can argue as in Case 1. Otherwise, by $(f3)$ (see Lemma 10), $|V(Q^{\downarrow}_*)|\geq2\delta-2\kappa+5$, whence
$$
|C|\geq |V^{\uparrow}|+|V(Q^{\downarrow}_*)|-2\geq(4\delta-2\kappa)+\delta-2>4\delta-2\kappa.
$$

\textbf{Case 3.} $f=2$ and $S\subseteq V^{\uparrow}.$

By Lemma 6, either $A^{\uparrow}-V^{\uparrow}$ is independent or $|V^{\uparrow}|\geq3\delta-5$. If $A^{\uparrow}-V^{\uparrow}$ is independent, then we can argue as in proof of Theorem 2 (Case 2.3). So, we can assume that $|V^{\uparrow}|\geq3\delta-5$. By Lemma 8, either $A^{\downarrow}-V^{\downarrow}=\emptyset$ or $|V^{\downarrow}|\geq2\delta-2\kappa+3$. If $A^{\downarrow}-V^{\downarrow}=\emptyset$, then
$$
|C|\geq |V^{\uparrow}|+|A^{\downarrow}|\geq(3\delta-5)+(\delta-\kappa+1)\geq4\delta-2\kappa.
$$

Let $|V^{\downarrow}|\geq2\delta-2\kappa+3$. Then 
$$
|C|\geq |V^{\uparrow}|+|V^{\downarrow}|-2\geq(3\delta-5)+(2\delta-2\kappa+3)-2\geq4\delta-2\kappa. \qquad   \Delta
$$

\noindent\textbf{Proof of Theorem 5}. If $\delta\leq2\kappa-3$, then we are done by Lemma 4. Let $\delta\geq2\kappa-2$ and let $Q^{\uparrow}_{1},...,Q^{\uparrow}_{m}$; $Q^{\downarrow}_{1},...,Q^{\downarrow}_{m}$ and $C$ be as defined in Definitions A-C. The existence of $Q^{\downarrow}_{1},...,Q^{\downarrow}_{m}$ follows from Lemma D. Assume w.l.o.g. that $F(Q^{\uparrow}_{i})=v_i$ and $L(Q^{\uparrow}_{i})=w_i$ $(i=1,...,m)$. By Lemma 6, either $A^{\uparrow}-V^{\uparrow}$ is independent or $|V^{\uparrow}|\geq3\delta-5$.\\ 

\textbf{Case 1.} $f\geq3$.

By Lemma 9, either $A^{\downarrow}-V^{\downarrow}$ is independent or $|V^{\downarrow}|\geq3\delta-3\kappa+f.$\\

\textbf{Case 1.1.} $A^{\uparrow}-V^{\uparrow}$ and $A^{\downarrow}-V^{\downarrow}$ both are independent sets.

If $A^{\uparrow}\subseteq V^{\uparrow}$ and $A^{\downarrow}\subseteq V^{\downarrow}$, then clearly $C$ is a dominating cycle. If $A^{\uparrow}\subseteq V^{\uparrow}$ and $A^{\downarrow}\not\subseteq V^{\downarrow}$, then by $(f2)$ (see Lemma 10), $|V^{\downarrow}|\geq2\delta-2\kappa+2f-m$, whence 
$$
c\geq|A^{\uparrow}|+|V^{\downarrow}|\geq(4\delta-2\kappa)+\delta-\kappa+f+m-3>4\delta-2\kappa.
$$

Next, if $A^{\uparrow}\not\subseteq V^{\uparrow}$ and $A^{\downarrow}\subseteq V^{\downarrow}$, then by $(f1)$ (see Lemma 10), $|V^{\uparrow}|\geq2\delta+m-2$, whence
$$
c\geq|V^{\uparrow}|+|A^{\downarrow}|\geq(4\delta-2\kappa)+\delta-\kappa+m>4\delta-2\kappa.
$$

Finally, if $A^{\uparrow}\not\subseteq V^{\uparrow}$ and $A^{\downarrow}\not\subseteq V^{\downarrow}$, then by $(f1)$ (see Lemma 10), $|V^{\uparrow}|\geq2\delta+m-2$ and by $(f2)$ (see Lemma 10), $|V^{\downarrow}|\geq2\delta-2\kappa+2f-m$. Hence,
$$
c\geq|V^{\uparrow}|+|V^{\downarrow}|-2m\geq(4\delta-2\kappa)+2f-2m-2>4\delta-2\kappa.
$$

\textbf{Case 1.2.} $A^{\uparrow}-V^{\uparrow}$ is independent and $|V^{\downarrow}|\geq3\delta-3\kappa+f$.

If $A^{\uparrow}\subseteq V^{\uparrow}$, then
$$
c\geq|A^{\uparrow}|+|V^{\downarrow}|\geq(4\delta-2\kappa)+2\delta-2\kappa+f-3\geq4\delta-2\kappa.
$$   

Otherwise, by $(f1)$ (see Lemma 10), $|V^{\uparrow}|\geq2\delta+m-2$ and 
$$
c\geq|V^{\uparrow}|+|V^{\downarrow}|-2m\geq(4\delta-2\kappa)+\delta-\kappa+f+m-2>4\delta-2\kappa.
$$

\textbf{Case 1.3.} $|V^{\uparrow}|\geq3\delta-5$ and $A^{\downarrow}-V^{\downarrow}$ is independent.

If $A^{\downarrow}\subseteq V^{\downarrow}$, then
$$
c\geq|V^{\uparrow}|+|A^{\downarrow}|\geq(4\delta-2\kappa)+2\delta-\kappa-3>4\delta-2\kappa.
$$

Otherwise, by $(f2)$ (see Lemma 10), $|V^{\downarrow}|\geq2\delta-2\kappa+2f-m$ and
$$
c\geq|V^{\uparrow}|+|V^{\downarrow}|-2m\geq(4\delta-2\kappa)+\delta+2f-3m-5\geq4\delta-2\kappa.
$$

\textbf{Case 1.4.} $|V^{\uparrow}|\geq3\delta-5$ and $|V^{\downarrow}|\geq3\delta-3\kappa+f$.

Clearly
$$
c\geq|V^{\uparrow}|+|V^{\downarrow}|-2m\geq(4\delta-2\kappa)+2\delta-\kappa-5\geq4\delta-2\kappa.
$$

\textbf{Case 2.} $f=2$ and $S\not\subseteq V^{\uparrow}$.

By Lemma 9, either $A^{\downarrow}-V(Q^{\downarrow}_*)$ is independent or $|V^{\downarrow}|\geq3\delta-3\kappa+3.$\\

\textbf{Case 2.1.} $A^{\uparrow}-V^{\uparrow}$ and $A^{\downarrow}-V(Q^{\downarrow}_*)$ both are independent sets.

If $A^{\uparrow}\subseteq V^{\uparrow}$ and $A^{\downarrow}\subseteq V^{\downarrow}$, then clearly $C$ is a dominating cycle. If $A^{\uparrow}\subseteq V^{\uparrow}$ and $A^{\downarrow}\not\subseteq V^{\downarrow}$, then by $(f3)$ (see Lemma 10), $|V^{\downarrow}|\geq2\delta-2\kappa+5$ and 
$$
c\geq|A^{\uparrow}|+|V^{\downarrow}|\geq(4\delta-2\kappa)+\delta-\kappa+2>4\delta-2\kappa.
$$

Further, if $A^{\uparrow}\not\subseteq V^{\uparrow}$ and $A^{\downarrow}\subseteq V^{\downarrow}$, then by $(f1)$ (see Lemma 10), $|V^{\uparrow}|\geq2\delta-1$, whence
$$
c\geq|V^{\uparrow}|+|A^{\downarrow}|\geq(4\delta-2\kappa)+\delta-\kappa+1>4\delta-2\kappa.
$$

Finally, if $A^{\uparrow}\not\subseteq V^{\uparrow}$ and $A^{\downarrow}\not\subseteq V^{\downarrow}$, then by $(f1)$ (see Lemma 10), $|V^{\uparrow}|\geq2\delta-1$ and by $(f3)$ (see Lemma 10), $|V^{\downarrow}|\geq2\delta-2\kappa+5$. Then
$$
c\geq|V^{\uparrow}|+|V^{\downarrow}|-2\geq(4\delta-2\kappa)+2>4\delta-2\kappa.
$$

\textbf{Case 2.2.} $A^{\uparrow}-V^{\uparrow}$ is independent and $|V^{\downarrow}|\geq3\delta-3\kappa+3$.

If  $A^{\uparrow}\subseteq V^{\uparrow}$, then
$$
c\geq|A^{\uparrow}|+|V^{\downarrow}|\geq(4\delta-2\kappa)+2\delta-2\kappa\geq4\delta-2\kappa.
$$   

Otherwise, by $(f1)$ (see Lemma 10), $|V^{\uparrow}|\geq2\delta-1$ and 
$$
c\geq|V^{\uparrow}|+|V^{\downarrow}|-2\geq(4\delta-2\kappa)+\delta-\kappa\geq4\delta-2\kappa.
$$

\textbf{Case 2.3.} $|V^{\uparrow}|\geq3\delta-5$ and $A^{\downarrow}-V(Q^{\downarrow}_*)$ is independent.

If $A^{\downarrow}\subseteq V(Q^{\downarrow}_*)$, then
$$
c\geq|V^{\uparrow}|+|A^{\downarrow}|\geq(4\delta-2\kappa)+2\delta-\kappa-3>4\delta-2\kappa.
$$

Otherwise, by $(f3)$ (see Lemma 10), $|V^{\downarrow}|\geq2\delta-2\kappa+5$ and
$$
c\geq|V^{\uparrow}|+|V^{\downarrow}|-2\geq(4\delta-2\kappa)+\delta-2>4\delta-2\kappa.
$$

\textbf{Case 2.4.} $|V^{\uparrow}|\geq3\delta-5$ and $|V^{\downarrow}|\geq3\delta-3\kappa+3$.

Clearly
$$
c\geq|V^{\uparrow}|+|V^{\downarrow}|-2\geq(4\delta-2\kappa)+2\delta-\kappa-4\geq4\delta-2\kappa.
$$

\textbf{Case 3.} $f=2$ and $S\subseteq V^{\uparrow}$.

By Lemma 6, either $A^{\uparrow}-V^{\uparrow}$ is independent or $|V^{\uparrow}|\geq3\delta-5$. By Lemma 9, $|V^{\downarrow}|\geq2\delta-2\kappa+3$.\\

\textbf{Case 3.1.} $A^{\uparrow}-V^{\uparrow}$ is independent. 

If $A^{\uparrow}\subseteq V^{\uparrow}$, then
$$
c\geq|A^{\uparrow}|+|V^{\downarrow}|\geq(4\delta-2\kappa)+\delta-\kappa\geq4\delta-2\kappa.
$$  

Otherwise, by $(f1)$ (see Lemma 10), $|V^{\uparrow}|\geq2\delta-1$ and
$$
c\geq|V^{\uparrow}|+|V^{\downarrow}|-2\geq4\delta-2\kappa.
$$

\textbf{Case 3.2.} $|V^{\uparrow}|\geq3\delta-5$.

Clearly
$$
c\geq|V^{\uparrow}|+|V^{\downarrow}|-2\geq(4\delta-2\kappa)+\delta-4\geq4\delta-2\kappa.  \qquad   \Delta
$$

\noindent\textbf{Proof of Theorem 1}. If $\delta\leq2\kappa-3$, then we are done by Lemma 4. Let $\delta\geq2\kappa-2$. Then the desired result follows from Theorems 2-5, immediately. \qquad $\Delta$

\end{document}